\documentclass{article}

\input xy
\xyoption {all}
\xyoption{v2}
\usepackage{amsmath, amsthm, amssymb}

\def\equals {\hspace{10pt}=\hspace{10pt}}

\def\Mor {\mathrm{Mor}}

\def\of #1{(#1)}

\def\del {\partial}

\newtheorem{definition}{Definition}

\newtheorem{proposition}{Proposition}
\newtheorem{theorem}{Theorem}
	
\newtheorem{observation}{Observation}

\def\proof {{Proof.}\hspace{7pt}}
\def\endofproof {\hfill $\Box$ \\}

\def\id {\mathrm{id}}

\def\gprod #1#2{{#2}{#1}}

\def\hprod #1#2{{#2}{#1}}

\def\act #1#2{{}^{#1}\!{#2}}
\def\leftact #1#2{#2}

\def\rightact #1#2{\act{#1}{#2}}

\begin{document}

\title{The inner automorphism 3-group of a strict 2-group}

\author{David Michael Roberts  and  Urs Schreiber}

\maketitle

\begin{abstract}
  Any group $G$ gives rise to a 2-group of inner automorphisms, $\mathrm{INN}(G)$.
  It is an old result by Segal that the nerve of this is the universal $G$-bundle. 
  We discuss that, similarly, for every 
  2-group $G_{(2)}$ there is a 3-group $\mathrm{INN}(G_{(2)})$  
  and a slightly smaller 3-group $\mathrm{INN}_0(G_{(2)})$
  of inner automorphisms. 
  We describe these for $G_{(2)}$ any strict 2-group, discuss how 
  $\mathrm{INN}_0(G_{(2)})$
  can be understood as arising from the mapping cone of the identity on
  $G_{(2)}$ and show that its underlying 2-groupoid structure fits 
  into a short exact sequence
  $$
    \xymatrix{
      G_{(2)}
      \ar[r]
      &
      \mathrm{INN}_0(G_{(2)})
      \ar[r]
      &
      \mathbf{B} G_{(2)}
    }
    \,.
  $$
  As a consequence, $\mathrm{INN}_0(G_{(2)})$ encodes the properties of the
  universal $G_{(2)}$ 2-bundle.
\end{abstract}

\tableofcontents

\section{Introduction}

The theory of groups and their principal fiber bundles generalizes
to that of categorical groups and their
categorical principal fiber bundles. In fact, using higher categories, one
has for each integer $n$ the notion of $n$-groups and their principal $n$-bundles.

The reader may have encountered principal 2-bundles mostly in the language
of (nonabelian) gerbes, which are to 2-bundles essentially like
sheaves are to ordinary bundles. The concept of a 2-bundle proper is described
in \cite{Bartels,Bakovic}.

These $n$-bundles are certainly interesting already in their own right.
One crucial motivation for considering them comes from the 
study of $n$-dimensional quantum field theory. In this case one is interested
in $n$-dimensional analogs of the concept of parallel transport in fiber 
bundles with connection \cite{BaezS,SWaldorfI,SWaldorfII}.

In that context a curious phenomenon occurs: 
whenever one investigates $n$-dimensional quantum field theory governed by
an $n$-group $G_{(n)}$, it turns out \cite{LieStruc}
that the situation is governed by an $(n+1)$-group associated to $G_{(n)}$.
In fact, it is appropriate to call this $(n+1)$-group  
$\mathrm{INN}_0(G_{(n)})$, because,
as the notation suggests, it is related to inner automorphisms
of the original $n$-group $G_{(n)}$. 

One of our aims here is to define
what inner automorphisms of a $2$-group are and to give a concise
definition as well as a detailed description of the 3-group 
$\mathrm{INN}_0(G_{(2)})$ for any strict 2-group $G_{(2)}$.
We then prove that $\mathrm{INN}_0(G_{(2)})$ has a couple of 
rather peculiar properties; it is \emph{contractible}
(equivalent, as a 2-group, to the trivial 2-group),
and fits into a short exact sequence
$$
  \xymatrix{
    G_{(2)}
    \ar@{^{(}->}[r]
    &
    \mathrm{INN}_0(G_{(2)})
    \ar@{->>}[r]
    &
    \mathbf{B} G_{(2)}
  }
$$
of 2-groupoids. To appreciate this result, it is helpful to first consider the
analogous statement for ordinary groups.

\paragraph{The statement for ordinary groups.}

For any ordinary group $G$, various constructions
of interest, like that of the universal $G$-bundle, are closely
related to a certain groupoid determined by $G$.

There are several different ways to think of this groupoid. 
The simplest way to describe its structure is to say that it is
the codiscrete groupoid
over the elements of $G$, namely the groupoid
whose objects are the elements of $G$ and which has precisely one
morphism from any element to any other.

The relevance of this groupoid is
better understood by thinking of it as the action groupoid
$
  G // G
$
of the action of $G$ on itself by left multiplication. As such, we may
write any of its morphisms as
$$
  \xymatrix{
    g \ar[r]^h & h g
  }
$$
for $g, h \in G$ and $hg$ being the product of $h$ and $g$ in $G$.

While this way of thinking about our groupoid already makes it 
more plausible that it is related to $G$-actions and hence 
possibly to $G$-bundles, one more property remains to be made manifest: 
there is also a monoidal structure on our groupoid. For any 
two morphisms,
$$
  \xymatrix{
    g_1 \ar[r]^{h_1} & h_1 g_1
  }
$$
and
$$
  \xymatrix{
    g_2 \ar[r]^{h_2} & h_2 g_2,
  }
$$
we can form the product morphism
$$
  \xymatrix{
    g_2 g_1 \ar[rr]^{h_2 \mathrm{Ad}_{g_2}(h_1) } && h_2 g_2 h_1 g_1,
  }
$$
and this assignment is functorial in both arguments. Moreover, to every morphism
$$
  \xymatrix{
    g \ar[r]^{h} & h g
  }
$$
there is a morphism 
$$
  \xymatrix{
    g^{-1} \ar[rr]^{\mathrm{Ad}_{g^{-1}}(h)^{-1}} && (h g)^{-1},
  }
$$
which is its inverse with respect to this product operation.

This makes $G // G$ a strict 2-group 
\cite{BaezLauda}. 
A helpful way to make the 2-group structure on
$G // G$ more manifest is to relate it to inner automorphisms of $G$.
To see this, consider another groupoid canonically associated to any group
$G$, namely the groupoid 
$$
  \mathbf{B} G = \lbrace \xymatrix{\bullet \ar[r]^g & \bullet}\; | g \in G\rbrace
$$
which has a single object $\bullet$, one morphism for each element of $G$ 
and where composition
of morphisms is just the product in the group.

Automorphisms
$$
  a : \mathbf{B} G \to \mathbf{B} G
$$
(i.e. invertible functors) of this groupoid are nothing but group automorphisms
of $G$. But now there are also isomorphisms between two such morphisms $a$ and $a'$, 
namely natural transformations:
$$
  \xymatrix{
    \mathbf{B} G
    \ar@/^1.8pc/[rr]^{a}_{\ }="s"
    \ar@/_1.8pc/[rr]_{a'}^{\ }="t"
    &&
    \mathbf{B} G
    \ar@{=>} "s"; "t"    
  }
  \,.
$$
This way for every ordinary group $G$ we have not just its ordinary group of
automorphisms, but actually a 2-group
$$
  \mathrm{AUT}(G) := \mathrm{Aut}_{\mathrm{Cat}}(\mathbf{B} G)
  \,.
$$
This is a groupoid, whose objects are group automorphisms of $G$.
The 2-group structure on this groupoid is manifest from the horizontal composition
of the natural transformations above.
Hence the ordinary automorphism group of $G$ is the group of objects of
$\mathrm{AUT}(G)$. 

By writing out the definition of a natural transformation,
one sees that there is a morphism between two objects in $\mathrm{AUT}(G)$ 
whenever the two underlying ordinary automorphisms of $G$ differ by conjugation 
with an element of $G$. It follows in particular that the \emph{inner}
automorphisms of $G$ correspond to those autofunctors of $\mathbf{B} G$
which are isomorphic to the identity:
$$
  \xymatrix{
    \mathbf{B} G
    \ar@/^1.8pc/[rr]^{\mathrm{Id}}_{\ }="s"
    \ar@/_1.8pc/[rr]_{\mathrm{Ad}_g}^{\ }="t"
    &&
    \mathbf{B} G
    \ar@{=>}^{\simeq}_g "s"; "t"    
  }
  \,.
$$

Therefore consider the groupoid
$
  \mathrm{INN}(G)
$:
its objects are pairs, consisting of an automorphisms together with
a transformation connecting it to the identity. A morphism from
$(g,\mathrm{Ad}_g)$ to $(gh,\mathrm{Ad}_{gh})$ is a commuting triangle
$$
  \xymatrix{
    &&
    \mathrm{Ad}_g
    \ar[dd]^{h}
    \\
    \mathrm{Id}_{\mathbf{B} G}
    \ar@/^1pc/[rru]^{g}_{\ }="s"
    \ar@/_1pc/[rrd]_{hg}^{\ }="t"
    \\
    &&
    \mathrm{Ad}_{hg}
    \ar@{=} "s"; "t" 
  }
  \,.
$$
This is again exactly the groupoid $G // G$ which we are discussing
$$
  \mathrm{INN}(G) = G//G
  \,.
$$
In this formulation the natural notion of composition of group automorphisms
nicely explains the monoidal structure on $G // G$.

Notice that $\mathrm{INN}(G)$ remembers the center of the group. We will
discuss that it sits inside an exact sequence
$$
  1 \to Z(G) \to \mathrm{INN}(G) \to \mathrm{AUT}(G) \to \mathrm{OUT}(G) \to 1
$$
of 2-groups, and that this is what generalizes to higher $n$.

If we think of the group $G$ just as a discrete category, whose objects are
the elements of $G$ and which has only identity morphisms,
then there is an obvious monomorphic functor
$$
  G \to \mathrm{INN}(G)
  \,.
$$
Moreover, there is an obvious epimorphic functor
$$
  \mathrm{INN}(G)
  \to \mathbf{B} G
$$
from our groupoid to the group $G$, but now with the latter regarded
as a category with a single object. This simply forgets the source and
target labels and recalls only the group element which is acting.

These two functors are such that the image of the former is precisely the
collection of morphisms which get sent to the identity morphism by the latter.
Therefore we say that we have a short exact sequence
\begin{eqnarray}
  \label{the short exact sequence of 1-groupoids}
  G \to \mathrm{INN}(G) \to \mathbf{B} G
\end{eqnarray}
of groupoids.

Notice that $G$ and $\mathrm{INN}(G)$ are groupoids which are also 2-groups
(the first one, being an ordinary group, is a degenerate case of a 2-group),
and that the morphism $G \to \mathrm{INN}(G)$ is also a morphism of
2-groups. But $\mathbf{B} G$ is in general just a groupoid 
without monoidal structure -- it has the structure of a 2-group
if and only if $G$ is abelian. 

Even though all this is rather elementary, the exact sequence 
(\ref{the short exact sequence of 1-groupoids}) is important.
We can apply the functor $|\cdot|$ to our sequence, which takes the nerve of a
category and then forms the geometric realization. Note that when $K$ is a 2-group, $|K|$ is a topological group.
Under this functor, (\ref{the short exact sequence of 1-groupoids})
becomes the universal $G$-bundle
$$
  G \to E G \to BG,
$$
even when $G$ is a topological or Lie group. The fact that \mbox{$BG \simeq |\mathbf{B} G|$} is the very definition of the classifying
space $BG$ of a group $G$ in \cite{Segal}. That \mbox{$EG \simeq |\mathrm{INN}(G)|$} 
is contractible follows from the existence of an equivalence of groupoids \mbox{$\mathrm{INN}(G) \stackrel{\sim}{\to} \ast$}. Finally, the inclusion \mbox{$G \to \mathrm{INN}(G)$} together with the 
monoidal structure on $\mathrm{INN}(G)$ gives the free $G$-action of $G$ on $EG$
whose quotient is exactly $B G$. The observation that $|\mathrm{INN}(G)|$ is a model
for $E G$ is originally due to Segal \cite{Segal}, who proved the remaining nontrivial
statement: that $|\mathrm{INN}(G)| \to B G$ is locally trivial when $G$ is a well-pointed group.\footnote{That is, the inclusion of the identity element is a closed cofibration.}

Our first main result is the higher categorical analogue of 
(\ref{the short exact sequence of 1-groupoids}),
obtained by starting with a strict 2-group $G_{(2)}$ in place of the
ordinary group $G$.

Here we do not consider geometric realizations of our categories
and 2-categories (for more on that see the closely related article
\cite{BaezStevenson} as well as \cite{RStevenson}) but instead focus on the existence of these sequences
of 1- and 2-groupoids. We comment on further aspects 
of the topic of universal $n$-bundles in \textsection\ref{universal n-bundles}. More
details will be given in \cite{Roberts}.

\paragraph{The formulation in terms of crossed modules.}

For many purposes, like doing explicit computations and for applying
the rich toolbox of simplicial methods, it is possible (and useful!) to express $n$-groups in terms of
$n$-term complexes of ordinary groups with extra structure on them. 
For instance strict 2-groups are well known to be equivalent to 
crossed modules of two ordinary groups: one describes the group of
objects, the other the group of morphisms of the 2-group.

This pattern continues, but there is a bifurcation of constructions, 
all of which are (homotopy) equivalent. 
Sufficiently strict 3-groups -- ``Gray groups'' -- are described 
by 2-crossed modules, which involve three ordinary groups forming a 
normal complex, and also by crossed squares, which look like crossed 
modules of crossed modules of groups. We will primarily use the former, 
and only mention crossed squares when we cannot avoid it.\footnote{As one goes to higher categorical dimensions (which we do not do here), there are multiple directions in which to extend the relevant diagrams, so there are a number of different models for $n$-groups. There is a sort of nonabelian Dold-Kan theorem, due to Carrasco and Cegarra \cite{Carrasco-Cegarra}, which can be used to characterise $n$-groups by $n$-term complexes of (possibly nonabelian) groups with the structure they call a hypercrossed complex.} The way we use these two models can be illustrated 
in one lower categorical dimension by comparing the \emph{map}
  $$
    \xymatrix{
      G \ar[r]^{\id} & G
    }
  $$
to the \emph{crossed module}
  $$
    \xymatrix{
      G \ar[r]^{\id} & G \ar[r]^<<<<<{\mathrm{Ad}} & \mathrm{Aut}(G)
    }
  $$
using that map. The crossed module can be thought of as the mapping cone ($=$homotopy quotient) of the identity map.
 
The translation between $n$-groups and their corresponding $n$-term
complexes of ordinary groups sheds light on both of these
points of view. The analogue of our statement about the 3-group $\mathrm{INN}_0(G_{(2)})$
is our second main result: the complex of groups describing
$\mathrm{INN}_0(G_{(2)})$ is the mapping cone of the identity on the
complex of groups describing $G_{(2)}$ itself.

This fact was anticipated from considerations in the theory of Lie $n$-algebras
\cite{LieStruc}, where the Lie $(n+1)$-algebra corresponding to a
Lie $(n+1)$-group $\mathrm{INN}_0(G_{(n)})$ has proven to be crucial 
for understanding connections with values in Lie $n$-algebras. There one
finds that inner derivation Lie $(n+1)$-algebras govern Lie $(n+1)$-algebras
of Chern-Simons type. The fact that $\mathrm{INN}_0(G_{(2)})$ arises from
a mapping cone of the identity is crucial in this context.

\paragraph{The plan of this article.}

The main content of this work is as follows -- first a concise and natural definition 
of inner automorphisms
of 2-groups, relating them to the full automorphism $(n+1)$-group
and to the categorical center. Then we apply this definition to the case
that $G_{(2)}$ is a strict 2-group and work out in full detail what
$\mathrm{INN}_0(G_{(2)})$ looks like, i.e. how the various composition
operations work, thus extracting 
its description in terms of complexes of ordinary groups. 
We state and prove the main properties of $\mathrm{INN}_0(G_{(2)})$.

The plan of our discussion is as follows.

\begin{itemize}

  \item
    In part \ref{n-groups in terms of groups} we recall the relation between
    2- and 3-groups and crossed modules of ordinary groups. This 
    serves to set up our convention for the precise choice of identification
    of 2-group morphisms with ordinary group elements.
    
  \item
    In part \ref{main results} we state our two main results.
    
  \item
    In part \ref{Inner automorphism n-groups} we define inner automorphism
    $n$-groups and prove some important general properties of them.
    
  \item
    In part \ref{The 3-group inng2} we apply our definition of inner automorphisms
    to an arbitrary strict 2-group $G_{(2)}$, to form the 3-group
    $\mathrm{INN}_0(G_{(2)})$. We then work out in detail the description of
    $\mathrm{INN}_0(G_{(2)})$ in terms of ordinary groups, spelling out
    the nature of the various composition and product operations.
    
  \item
    In part \ref{properties} we state and prove the main properties of
    $\mathrm{INN}_0(G_{(2)})$, including our two main results.
    
  \item
    In part \ref{universal n-bundles} we close by indicating in more detail
    how inner automorphism $(n+1)$-groups play the role of 
    universal $n$-bundles. We also relate our construction here to analogous
    simplicial constructions.

\end{itemize}

\medskip

We are grateful to Jim Stasheff for helpful discussions and for
emphasizing the importance of the mapping cone construction in the
present context. 
We profited from general discussion with Danny Stevenson 
and thank him for his help on the references to Segal's work.
We also thank 
Christoph Schweigert and
Zoran {\v S}koda for helpful comments on the manuscript,
Tim Porter for reminding us of Norrie's work.
We are grateful to Todd Trimble for discussion about the simplicial
aspects of our construction and its relation to d{\'e}calage,
and to an anonymous referee for making the remark which 
appears in \ref{Universal simplicial bundles}.

DMR is supported by an Australian Postgraduate Award.

\section{$n$-Groups in terms of groups}
\label{n-groups in terms of groups}

Sufficiently strict $n$-groups are equivalent to certain structures
-- crossed modules and generalizations theoreof -- involving
just collections of ordinary groups with certain structure on them.

\subsection{Conventions for strict 2-groups and crossed modules}

An ordinary group $G$ may be regarded as a one object category. If we
regard $G$ as such a category, we write $\mathbf{B} G$ in order to
emphasize that we are thinking of the monoidal 0-category $G$ as
a one object 1-category.

This way we obtain a notion of $n$-groups from any notion of $n$-categories:
an $n$-group $G_{(n)}$ is a monoidal $(n-1)$-category such that when
regarded as a one-object $(n)$-category $\mathbf{B} G_{(n)}$ it becomes
a one-object $n$-groupoid. An $n$-groupoid is an $n$-category all whose
$k$-morphisms are equivalences, for all $1 \leq k \leq n$.

Here we shall be concerned with strict 2-groups and with 3-groups which are
Gray-categories. A strict 2-group $G_{(2)}$ is one such that $\mathbf{B} G_{(2)}$
is a strict one-object 2-groupoid. A Gray-groupoid is a 3-groupoid which 
is strict except for the exchange law of 2-morphisms.

The standard reference for 2-groups is \cite{BaezLauda}. A discussion of 
Gray-groupoids useful for our context is in \cite{KampsPorter}.  We come
to Gray-groups in \ref{3-groups and 2-crossed modules}.

It is well known that strict 
2-groups are equivalent to crossed modules of ordinary groups.
This was first established in \cite{BrownSpencer}. The relation to 
category objects in groups was also discussed in \cite{Loday}.
The notion of a crossed module is originally due to \cite{Whitehead}.

\begin{definition}
  A {\bf crossed module of groups} is a diagram
  $$
    \xymatrix{
      H \ar[r]^t & G \ar[r]^<<<<<\alpha & \mathrm{Aut}(H)
    }
  $$
  in $\mathrm{Grp}$ such that
  $$
    \raisebox{20pt}{
    \xymatrix{
      H 
       \ar[dr]_t 
       \ar[rr]^{\mathrm{Ad}}
       && 
       \mathrm{Aut}(H)
      \\
      & G \ar[ur]_\alpha 
    }
    }
  $$
  and
  $$
    \raisebox{36pt}{
    \xymatrix{
      G \times H
      \ar[rr]^{\mathrm{Id} \times t}
      \ar[dd]_{\alpha}
      &&
      G \times G
      \ar[dd]^{\mathrm{Ad}}
      \\
      \\
      H
      \ar[rr]^t
      &&
      G
    }
    }
    \,.
  $$
\end{definition}

\begin{definition}
  A {\bf strict 2-group} $G_{(2)}$ is any of the following equivalent
  entities
  \begin{itemize}
    \item a group object in $\mathrm{Cat}$
    \item a category object in $\mathrm{Grp}$
    \item a strict 2-groupoid with a single object
  \end{itemize}  
\end{definition}
A detailed discussion can be found in \cite{BaezLauda}.

One identifies 
\begin{itemize}
  \item
    $G$ is the group of objects of $G_{(2)}$.
  \item
    $H$ is the group of morphism of $G_{(2)}$ starting at the identity object. 
  \item
    $t : H \to G$ is the target homomorphism so that
    $h : \mathrm{Id} \to t(h)$ for all $h \in H$.
  \item
    $\alpha : G \to \mathrm{Aut}(H)$ is conjugation with identity morphisms:
    $$ 
      \mathrm{Ad}_{\mathrm{Id}_g}(\xymatrix{\mathrm{Id} \ar[r]^h & t(h)})
      =
      \xymatrix{\mathrm{Id} \ar[rr]^<<<<<<<<<<<{\alpha(h)(h)} && t(\alpha(g)(h))}
    $$
    for all $g \in G$, $h \in H$.
\end{itemize}

We often abbreviate
$$
  \act{g}{h} := \alpha(g)(h) 
  \,.
$$

Beyond that there are $2 \times 2$ choices to be made when identifying a 
strict 2-group $G_{(2)}$ with a crossed module 
of groups.

The first choice to be made is
in which order to multiply elements in $G$. For 
$\xymatrix{
  \bullet \ar[r]^{g_1} & \bullet
 }$ and 
$\xymatrix{
  \bullet \ar[r]^{g_2} & \bullet
 }$
two morphisms in $\mathbf{B} G_{(2)}$, we can either set

\begin{eqnarray*}
  \xymatrix{
    \bullet \ar[r]^{g_1} & \bullet \ar[r]^{g_{2}} & \bullet
    \;
    :=
    \;
    \bullet \ar[rr]^{g_1 g_2} && \bullet
  }
  &&
  \mbox{(F)}
\end{eqnarray*}
or
\begin{eqnarray*}
  \xymatrix{
    \bullet \ar[r]^{g_1} & \bullet \ar[r]^{g_{2}} & \bullet
    \;
    :=
    \;
    \bullet \ar[rr]^{g_2 g_1} && \bullet
  }
  &&
  \mbox{(B)}  
  \,.
\end{eqnarray*}

The other choice to be made is how to describe arbitrary morphisms by
an element in the semidirect product group $G \ltimes H$: 
every morphism of $G_{(2)}$ may be written 
as the product of one starting at the identity object with an identity morphism
on some object. The choice of ordering here yields either

\begin{eqnarray*}
  \xymatrix{
    \bullet
    \ar@/^1.8pc/[rr]^{g}_{\ }="s"
    \ar@/_1.8pc/[rr]_{}^{\ }="t"
    &&
    \bullet
    \ar@{=>}^h "s"; "t"
  }
  \;
  :=
  \;
  \xymatrix{
    \bullet
    \ar@/^1.8pc/[rr]^{\mathrm{Id}}_{\ }="s"
    \ar@/_1.8pc/[rr]_{}^{\ }="t"
    &&
    \bullet
    \ar[rr]^{g}
    &&
    \bullet
    \ar@{=>}^h "s"; "t"
  }
  &&
  \mbox{(R)}
\end{eqnarray*}
or
\begin{eqnarray*}
  \xymatrix{
    \bullet
    \ar@/^1.8pc/[rr]^{g}_{\ }="s"
    \ar@/_1.8pc/[rr]_{}^{\ }="t"
    &&
    \bullet
    \ar@{=>}^h "s"; "t"
  }
  \;
  :=
  \;
  \xymatrix{
    \bullet
    \ar[rr]^{g}
    &&
    \bullet
    \ar@/^1.8pc/[rr]^{\mathrm{Id}}_{\ }="s"
    \ar@/_1.8pc/[rr]_{}^{\ }="t"
    &&
    \bullet
    \ar@{=>}^h "s"; "t"
  }
  &&
  \mbox{(L)}
\end{eqnarray*}

Here we choose the convention 
$$LB\,.$$
This implies 
$$
  \xymatrix{
    \bullet
    \ar@/^1.8pc/[rr]^{g}_{\ }="s"
    \ar@/_1.8pc/[rr]_{g' = t(h)g}^{\ }="t"
    &&
    \bullet
    \ar@{=>}^h "s"; "t"
  }
$$
for all $g \in G$, $h \in H$, 
as well as the following two equations for
horizontal and vertical composition in $\mathbf{B} G_{(2)}$, expressed in
terms of operations in the crossed module

$$
  \xymatrix{
    \bullet
    \ar@/^1.8pc/[rr]^{g_1}_{\ }="s1"
    \ar@/_1.8pc/[rr]_{}^{\ }="t1"
    &&
    \bullet
    \ar@/^1.8pc/[rr]^{g_2}_{\ }="s2"
    \ar@/_1.8pc/[rr]_{}^{\ }="t2"
    &&
    \bullet
    \ar@{=>}^{h_1} "s1"; "t1"
    \ar@{=>}^{h_2} "s2"; "t2"
  }
  =
  \xymatrix{
    \bullet
    \ar@/^1.8pc/[rr]^{\gprod{g_1}{g_2}}_{\ }="s"
    \ar@/_1.8pc/[rr]_{}^{\ }="t"
    &&
    \bullet
    \ar@{=>}|{ \hprod{\rightact{g_2}{h_1}}{h_2}} "s"; "t"
  }  
$$
and
$$
  \xymatrix{
    \bullet
    \ar@/^2.3pc/[rr]^{g_1}_{\ }="s"
    \ar[rr]|{g_2}="c"
    \ar@/_2.3pc/[rr]_{g_3}^{\ }="t"
    &&
    \bullet
    \ar@{=>}^{h_1} "s"; "c"
    \ar@{=>}^{h_2} "c"; "t"
  }
  =
  \xymatrix{
    \bullet
    \ar@/^1.8pc/[rr]^{g_1}_{\ }="s"
    \ar@/_1.8pc/[rr]_{}^{\ }="t"
    &&
    \bullet
    \ar@{=>}^{\hprod{h_1}{h_2}} "s"; "t"
  } 
  \,.
$$

\subsection{3-Groups and 2-crossed modules}
\label{3-groups and 2-crossed modules}

As we are considering strict models in this paper, 
we will assume that all 3-groups are as strict as possible. 
This means they will be one-object 
Gray-categories, or \emph{Gray-monoids} \cite{Day-Street}. 
A Gray-monoid is a (strict) 2-category $\cal{M}$ such that the product 2-functor
$$\cal{M} \otimes \cal{M} \to \cal{M}$$
uses the Gray tensor product \cite{Gray,GordonPowerStreet}, 
not the usual Cartesian product of 2-categories. Thus non-identity coherence morphisms only appear when we use the monoidal structure on $\cal{M}$. 
So from now on a ``3-group'' $G_{(3)}$ will mean a 3-group such that
regarded as a one-object 3-groupoid $\mathbf{B} G_{(3)}$ it is a one-object Gray-groupoid.

Just as a 2-group gives rise to a crossed module, a 3-group gives rise to a 2-crossed module. 
Roughly, this is a complex of groups

$$
  L \to M \to N,
$$
and a function
\begin{eqnarray}
  \label{Peiffer lift}
      M \times M \to L           
\end{eqnarray}
such that $L \to M$ is a crossed module, and (\ref{Peiffer lift}) 
measures the failure of $M \to N$ to be a crossed module. An example is
when $L = 1$, and then we have a crossed module. 

The relation between 3-groups and 2-crossed modules was described in 
\cite{KampsPorter}.
The precise definition of a 2-crossed module is as follows, 
see also \cite{Conduche}.

\begin{definition}\label{2xmod definition}
A {\bf 2-crossed module} is a normal complex of length 2

$$
\xymatrix{
L \ar[r]^{\del_2} & M \ar[r]^{\del_1} & N
}
$$
of $N$-groups ($N$ acting on itself by conjugation) and an $N$-equivariant function
$$\{\cdot,\cdot\}: M \times M \to L\,,$$
called a Peiffer lifting, satisfying these conditions:

\begin{enumerate}
\item
$\del_2\{m,m'\} = (mm'm^{-1})(\act{\del_1m}{m'})^{-1},$
\item
$\{\del_2l,\del_2l'\} = [l,l'] := ll'l^{-1}l'^{-1}$
\item
\begin{enumerate}
   \item
   $\{m,m'm''\} = \{m,m'\}\act{mm'm^{-1}}{\{m,m''\}}$
   \item
   $\{mm',m''\} = \{m,m'm''m'^{-1}\}\act{\del_1m}{\{m',m''\}},$
\end{enumerate}
\item
$\{m,\del_2l\} = (\act{m}{l})(\act{\del_1m}{l})^{-1}$
\item
$\act{n}{\{m,m'\}} = \{\act{n}{m},\act{n}{m'}\},$
\end{enumerate}
 where $l,l'\in L,\ m,m',m''\in M$ and $n\in N$. 
\end{definition}

Here ${}^m l$ denotes the action
\begin{eqnarray}
   \label{M acting on L}
   M \times L &\to& L\nonumber\\
   (m,l) &\mapsto& \act{m}{l} := l\{\del_2l^{-1},m\}
   \,.
\end{eqnarray}

A \emph{normal} 
complex is one in which $im\ \partial$ is normal in $ker\ \partial$ 
for all differentials.

It follows from these conditions that $\del_2:L \to M$ is a crossed module with the action (\ref{M acting on L}).

To get from a 3-group $G_{(3)}$ to a 2-crossed module \cite{KampsPorter}, we emulate the construction of a crossed module from a 2-group: one identitfies 
\begin{itemize}
  \item
    $N$ is the group of objects of $G_{(3)}$.
  \item
    $M$ is the group of 1-morphisms of $G_{(3)}$ starting at the identity object. 
  \item
    $L$ is the group of 2-morphisms starting at the identity 1-arrow of the identity object
  \item
    $\del_1 : M \to N$ is the target homomorphism such that
    $m : \mathrm{Id} \to \del_1(m)$ for all $m \in M$.
  \item
    $\del_2 : L \to M$ is the target homomorphism such that
    $l : \mathrm{Id}_{\mathrm{Id}} \to \del_2(l)$ for all $l \in L$.
  \item
    The various actions arise by whiskering, analogously to the case of a 2-group.
\end{itemize}

We will not go into the proof that this gives rise 
to a 2-crossed module for all 3-groups, 
but only in the case we are considering. One reason to consider 
2-crossed modules is that the homotopy groups of $G_{(3)}$ 
can be calculated as the homology of the sequence underlying the 2-crossed module.

\subsection{Mapping cones of crossed modules}
\label{Mapping cones of crossed modules}

Another notion related to 3-groups 
\cite{ArvasiUlualan}
is crossed modules internal to 
crossed modules (more technically known as crossed squares, \cite{Guin-Walery--Loday},\cite{Loday}). 
More generally, consider a map $\phi$ of crossed modules:

\begin{definition}[nonabelian mapping cone \cite{Loday}]
  \label{nonabelian mapping cone}
  For
  $$
    \xymatrix{
      H_2 
       \ar[rr]^{\phi_H}
       \ar[dd]^{t_2}
       && 
       H_1
       \ar[dd]^{t_1}
      \\
      \\
      G_2 \ar[rr]^{\phi_G}&& G_1 
    }
  $$
  a 2-term complex of crossed modules ($t_i : H_i \to G_i$),
  we say its mapping cone is the complex of groups
  \begin{eqnarray}
   \label{mapping cone}
    \xymatrix{
      H_2 
        \ar[rr]^<<<<<<<<<<{\partial_2} 
      && 
      G_2 \ltimes H_1   
      \ar[rr]^{\partial_1}
      &&
      G_1
    }
    \,,
\end{eqnarray}
  where
  $$
    \partial_1 : (g_2,h_1) \mapsto t_1(h_1)\phi_G(g_2)
  $$
  and
  $$
    \partial_2 : h_2 \mapsto (t_2(h_2),\phi_H(h_2)^{-1})
    \,.
  $$
\end{definition}
Here $G_2$ acts on $H_1$ by way of the morphism $\phi_G : G_2 \to G_1$.

When no structure is imposed on $\phi$, (\ref{mapping cone}) is merely a complex. However, if $\phi$ is a crossed square, the mapping cone is a 2-crossed module (originally shown in \cite{Conduche letter}, but see \cite{Conduche}). 
We will not need to define crossed squares here 
but just note they come equipped with a map
$$
  h : G_2 \times H_1 \to H_2
$$
satisfying conditions similar to the Peiffer lifting.
The details can be found in \cite{Conduche}, which we recall in our appendix. 
The only crossed square we will see in this paper is the identity map on a crossed module
$$
  \xymatrix{
   H 
       \ar[rr]^{\mathrm{id}}
       \ar[dd]_{t}
       && 
       H
       \ar[dd]^{t}
      \\
      \\
      G \ar[rr]^{\mathrm{id}}&& G 
    }
  $$
with the structure map
\begin{eqnarray}
  h: G \times H &\to& H\\
       (g,h) &\mapsto & h\act{g}{h^{-1}}.
\end{eqnarray}

This concept of ``crossed modules of crossed modules'' is explored in Norrie's 
thesis \cite{Norrie} on `actors' of crossed modules, with a focus on categorifying group theory, 
rather than geometry. Automorphisms of crossed modules of groups and groupoids have
been discussed in \cite{BrownGilbert} and \cite{BrownIcen}.

\section{Main results}
\label{main results}

\subsection{The exact sequence 
$G_{(2)} \to \mathrm{INN}_0(G_{(2)}) \to \mathbf{B} G_{(2)}$}
\label{the exact sequence}

We describe the 3-group 
$\mathrm{INN}_0(G_{(2)})$ for $G_{(2)}$ any strict 
2-group, and show that it plays the role of 
the universal principal $G_{(2)}$-bundle
in that
\begin{itemize}
  \item
    $\mathrm{INN}_0(G_{(2)})$ is equivalent to the trivial 3-group 
    (hence ``contractible'').
  \item
    $\mathrm{INN}_0(G_{(2)})$ fits into the short exact sequence
  $$
    \xymatrix{
      G_{(2)}
      \ar@{^{(}->}[r]
      &
      \mathrm{INN}_0(G_{(2)})
      \ar@{->>}[r]
      &
      \mathbf{B} G_{(2)}
    }
  $$
  of strict 2-groupoids.
\end{itemize}

\subsection{$\mathrm{INN}_0(G_{(2)})$ from a mapping cone}

We show that the 3-group $\mathrm{INN}_0(G_{(2)})$ comes
from a 2-crossed module 
$$
 \xymatrix{
  H \ar[r] 
  &
  G \ltimes H
  \ar[r]
  &
  G
 }
$$
which is the mapping cone of  
$$
  \raisebox{20pt}{
  \xymatrix{
    H 
      \ar[rr]^{\mathrm{Id}}
      \ar[dd]_{t}
    &&
    H 
    \ar[dd]^t
    \\
    \\
    G \ar[rr]^{\mathrm{Id}}
    &&
    G
  }}\;,
$$
the identity map of the crossed module
$(t : H \to G)$ which determines $G_{(2)}$.

Notice that this harmonizes with the analogous 
result for Lie 2-algebras
discussed in \cite{LieStruc}.

\section{Inner automorphism $(n+1)$-groups}
\label{Inner automorphism n-groups}

An automorphism of an $n$-group $G_{(n)}$ is 
simply an automorphism of the $n$-category $\mathbf{B} G_{(n)}$.
We want to say that such an automorphism $q$ is \emph{inner} if it is
equivalent to the identity automorphism
$$
  \xymatrix{
    \mathbf{B} G_{(n)}
    \ar@/^1.8pc/[rr]^{\mathrm{Id}}_{\ }="s"
    \ar@/_1.8pc/[rr]_{q}^{\ }="t"
    &&
    \mathbf{B} G_{(n)}
    \ar@{=>}^\sim "s"; "t"
    \,.
  }
$$
Notice that we really do mean automorphisms here, and not auto-equivalences: we require
an autmorphism to be an endo-$n$-functor with a \emph{strict} inverse. 

Automorphisms (of crossed modules) connected to the identity appear in definition 2.3 of
\cite{BrownIcen} under the name ``free derivations''. Since the naturality diagram for
the transformation connecting an automorphism to the identity implies that this automorphism
arises from conjugations, we think of them as \emph{inner automorphisms} here and 
reserve the term ``(inner) derivations'' for the image of these automorphisms as one passes from
Lie $n$-groups to Lie $n$-algebras \cite{LieStruc}.

A useful way to think of the $n$-groupoid of inner automorphisms
is in terms of what we call ``tangent categories'', a slight variation
of the concept of comma categories.

Tangent categories in general happen to live in interesting exact sequences.
In order to be able to talk about these, we first quickly set up a
our definitions for exact sequences of strict 2-groupoids.

Remember that we work entirely within
the Gray-category whose objects are strict 2-groupoids,
whose morphisms are strict 2-functors, whose 2-morphisms are
pseudonatural transformations and whose 3-morphisms are 
modifications of these.

\subsection{Exact sequences of strict $2$-groupoids}

Inner automorphism $n$-groups turn out to live in interesting
\emph{exact sequences of $(n+1)$-groups}.
Therefore we want to talk about generalizations of exact sequences of groups
to the world of $n$-groupoids. Since for our purposes here only
strict 2-groupoids matter, we shall be content with just using
a definition applicable to that case.

\begin{definition}[exact sequence of strict 2-groupoids]
  A collection of composable morphisms
  $$
    \xymatrix{
      C_0 
      \ar[r]^{f_1}
      &
      C_1 
      \ar[r]^{f_2}
      &
      \cdots
      \ar[r]^{f_n}
      &
      C_n
    }
  $$
  of strict 2-groupoids $C_i$
  is called an exact sequence
  if, as ordinary maps between spaces of 2-morphisms, 
  \begin{itemize}
    \item
       $f_1$ is injective
    \item
       $f_n$ is surjective
    \item
      the image of $f_i$ is the preimage under $f_{i+1}$ of the 
      collection of all identity 2-morphisms on identity 1-morphisms
      in $\mathrm{Mor}_2(C_{i+1})$, for all $1 \leq i < n$.
  \end{itemize}
\end{definition}

In order to make this harmonize with our distinction between 
$n$-groups $G_{(n)}$ and the corresponding 1-object $n$-groupoids
$\mathbf{B} G_{(n)}$ we add to that

\begin{definition}[exact sequences of strict 2-groups]
  A collection of composable morphisms
  $$
    \xymatrix{
      G_0 
      \ar[r]^{f_1}
      &
      G_1 
      \ar[r]^{f_2}
      &
      \cdots
      \ar[r]^{f_n}
      &
      G_n
    }
  $$
  of strict 2-groups is called an exact sequence if the corresponding
  chain
  $$
    \xymatrix{
      \mathbf{B} G_0 
      \ar[r]^{\mathbf{B} f_1}
      &
      \mathbf{B} G_1 
      \ar[r]^{\mathbf{B} f_2}
      &
      \cdots
      \ar[r]^{\mathbf{B} f_n}
      &
      \mathbf{B} G_n
    }
  $$
  is an exact sequence of strict 2-groupoids.
\end{definition}

\paragraph{Remark.}
Ordinary exact sequences of groups thus precisely correspond
to exact sequences of strict 2-groups all whose morphisms are
identities.

\subsection{Tangent $2$-categories}

We present a simple but useful way describe $2$-categories
of morphisms with coinciding source. We find it helpful to
refer to this construction as \emph{tangent categories}
for reasons to become clear. 
It is not hard to see that 
this is the globular analog of the simplicial construction 
known as \emph{d{\'e}calage}, as will be discussed more in
\ref{relation to simplicial bundles}. In the context of higher
categories it has been considered (over single objects)
in section 3.2 of \cite{KampsPorter}.

\begin{definition}
  Denote by
  $$
    \mathrm{pt} := \{\bullet\}
  $$
  the strict 2-category 
  with a single object and no nontrivial morphisms
  and by 
  $$
    I := \{ \xymatrix{\bullet \ar[r]^\sim & \circ }\}
    \,,
  $$
  the strict 2-category 
  consisting of two objects connected by a 1-isomorphism.
\end{definition}
Of course $I$ is equivalent to $\mathrm{pt}$ -- but not isomorphic.
We fix one injection
$$
  i : \xymatrix{
    \mathrm{pt} \ar@{^{(}->}[r] & I
  }
$$
$$
  i : \bullet \mapsto \bullet
$$
once and for all.

It is useful to think of morphisms 
$$
  \mathbf{f} : I \to C
$$
from $I$ to some codomain $C$ as labeled by the corresponding image
of the ordinary point
$$
  \raisebox{40pt}{
  \xymatrix{
    \mathrm{pt}
    \ar@{^{(}->}[dd]
    \ar[rr]^f
    &&
    C
    \ar[dd]^{=}
    \\
    \\
    I
    \ar[rr]^{\mathbf{f}}
    &&
    \mathrm{C}
  }
  }
  \,.
$$

\begin{definition}[tangent $2$-bundle]\label{tangent 2-bundle}
  Given any strict $2$-category $C$, we define its tangent 2-bundle
  $$
    T C \subset \mathrm{Hom}_{2\mathrm{Cat}}(I,C)
  $$
  to be that sub $2$-category of morphisms from $I$ into $C$
  which collapses to a 0-category when pulled back along the fixed inclusion
$
  i : \xymatrix{
    \mathrm{pt} \ar@{^{(}->}[r] & I
  }
$:
the morphisms $h$ in T C are all those for which
$$
  \raisebox{30pt}{
  \xymatrix{
    \mathrm{pt}
    \ar@{^{(}->}[dd]
    \\
    \\
    I
    \ar@/^1.8pc/[rr]^{\mathbf{f}}_{\ }="s"
    \ar@/_1.8pc/[rr]_{\mathbf{f}'}^{}="t"
    &&
    C
    \ar@{=>}^h "s"; "t"
  }
  }
  \hspace{7pt}
  =
  \hspace{7pt}
  \raisebox{30pt}{
  \xymatrix{
    \mathrm{pt}
    \ar@{^{(}->}[dd]
    \\
    \\
    I
    \ar[rr]^{\mathbf{f}}
    &&
    C
  }
  }
  \,.
$$
The tangent 2-bundle is a disjoint union
$$
  TC 
    = 
    \bigoplus_{x \in \mathrm{Obj}(C)}
    T_x C
$$
of tangent $2$-categories at each object $x$ of $C$. In this way
it is a $2$-bundle
$$
  \xymatrix{
    p : TC \ar[r] & \mathrm{Obj}(C)
  }
$$
over the space of objects of $C$.  
\end{definition}
As befits a tangent bundle, the tangent $2$-bundle has a canonical section
$$
  e_{\mathrm{Id}} : \mathrm{Obj}(C) \to TC 
$$
which sends every object of $C$ to the Identity morphism on it.

\paragraph{Remark.} The groupoid $I$ plays the role of the interval in topology
and underlies the homotopy theory for groupoids, as described 
in \cite{KampsPorterBook}. The terminology ``tangent category'' finds further
justification when the discussion here is done for smooth $n$-groups which are 
then sent to the corresponding Lie $n$-algebras: indeed, as indicated in figure
3 of \cite{LieStruc}, one finds a close relation between maps from the interval $I$,
inner automorphisms, the notion of universal $n$-bundles and tangency relations.

\paragraph{Example (slice categories).}
For $C$ any 1-groupoid, i.e. a strict 2-groupoid with only
identity 2-morphisms, its tangent 1-category is the comma category
$$
  TC
  =
  ( (\mathrm{Obj}(C) \hookrightarrow C) \downarrow \mathrm{Id}_{C})
  \,.
$$ 
This is the disjoint union of all 
co-over categories on all objects of $C$
$$
  TC = \bigoplus_{a \in \mathrm{Obj}(C)} (a \downarrow C)
$$
Objects of $TC$ are morphisms
$f : a \to b$ in $C$, 
and morphisms $\xymatrix{f \ar[r]^h & f'}$ in 
$TC$ are commuting triangles
$$
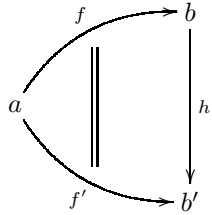
\raisebox{35pt}{
\xymatrix{
  &&
  b
  \ar[dd]^{h}
  \\
  a 
  \ar@/^1pc/[urr]^{f}_{\ }="s" 
  \ar@/_1pc/[drr]^{\ }="t"_{f'}
  &&
  \\
  &&
  b' 
  \ar@{=} "s";"t"
}
}
$$
in $C$.

\paragraph{Example (strict tangent 2-groupoids).}
The example which we are mainly interested in is that
where $C$ is a strict 2-groupoid. For $a$ any object
in $C$, an object of $T_a C$ is a morphism
$$
  \xymatrix{
    a \ar[r]^q & b
    \,.
  }
$$
A 1-morphism in $T_a C$ is a filled triangle
$$
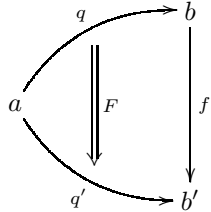
\raisebox{35pt}{
\xymatrix{
  &&
  b
  \ar[dd]^{f}
  \\
  a 
  \ar@/^1pc/[urr]^{q}_{\ }="s" 
  \ar@/_1pc/[drr]^{\ }="t"_{q'}
  &&
  \\
  &&
  b' 
  \ar@{=>}^F "s";"t"
}
}
$$
in $C$. Finally, a 2-morphism in $T_a C$ looks like
$$
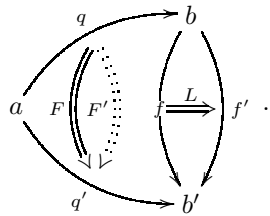
\raisebox{36pt}{
\xymatrix{
&&
  b 
  \ar @/^1pc/ [dd]^{f'}_{}="t_2"
  \ar @/_1pc/ [dd]|{f}^{}="s_2"
\\
  a 
 \ar @/^1pc/ [urr]^{q}_{\ }="s" \ar @/_1pc/ [drr]^{\ }="t"_{q'}&&
\\
&&
  b' 
\ar @{=>} @/_.7pc/ "s";"t"_{F} |>\hole
\ar @{:>} @/^.7pc/ "s";"t"_{F'} |>\hole
\ar @2 "s_2";"t_2"^{L}
}
}
\,.
$$
The composition of these 2-morphisms is the obvious one. 
We give a detailed description for the case the
$C = \mathbf{B} G_{(2)}$ in \ref{The 3-group inng2}.

\begin{proposition}
 \label{exact sequence}
 For any strict $2$-category $C$, its tangent $2$-bundle $TC$ fits 
 into an exact sequence
 $$
  \xymatrix{
    \mathrm{Mor}(C)
    \ar@{^{(}->}[r]
    &
    TC
    \ar@{->>}[r]
    &
    C
  }
 $$
 of strict 2-categories.
\end{proposition}
Here $\mathrm{Mor}(C) := \mathrm{Disc}(\mathrm{Mor}(C))$ is the 1-category
of morphisms of $C$, regarded as a strict 2-category with only identity 2-morphisms.

\proof
The strict inclusion 2-functor on the left is
$$
  \left(
  
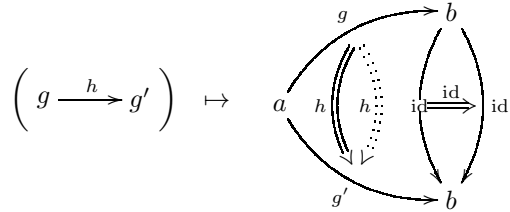
\xymatrix{
    g
    \ar[r]^h
    &
    g'
  }
  \right)
  \hspace{7pt}
    \mapsto
  \hspace{7pt}
\raisebox{36pt}{
\xymatrix{
&&
  b 
  \ar @/^1pc/ [dd]^{\id}_{}="t_2"\ar @/_1pc/ [dd]|{\id}^{}="s_2"
\\
  a 
 \ar @/^1pc/ [urr]^{g}_{\ }="s" \ar @/_1pc/ [drr]^{\ }="t"_{g'}&&
\\
&&
  b 
\ar @{=>} @/_.7pc/ "s";"t"_{h} |>\hole
\ar @{:>} @/^.7pc/ "s";"t"_{h} |>\hole
\ar @2 "s_2";"t_2"^{\id}
}
}
$$
for $g,g' : a \to b$ any two parallel morphisms in $C$ and $h$ any 2-morphism between them.

The strict surjection 2-functor on the right is
$$
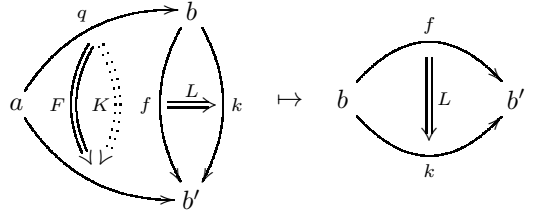
\raisebox{36pt}{
\xymatrix{
&&b \ar @/^1pc/ [dd]^{k}_{}="t_2"\ar @/_1pc/ [dd]_{f}^{}="s_2"
\\
a \ar @/^1pc/ [urr]^{q}_{\ }="s" \ar @/_1pc/ [drr]^{\ }="t"_{}&&
\\
&&b' 
\ar @{=>} @/_.7pc/ "s";"t"_{F} |>\hole
\ar @{:>} @/^.7pc/ "s";"t"_{K} |>\hole
\ar @2 "s_2";"t_2"^{L}
}
}
\hspace{7pt}
  \mapsto
\hspace{7pt}
\xymatrix{
  b
  \ar@/^1.8pc/[rr]^{f}_{\ }="s"
  \ar@/_1.8pc/[rr]_{k}^{\ }="t"
  &&
  b'
  \ar@{=>}^L "s"; "t"  
}
\,.
$$
The image of the injection is precisely the preimage under the surjection of
the identity 2-morphism on the identity 1-morphisms . This
means the sequence is exact.
\endofproof

\subsection{Inner automorphisms}

Often, for $G$ any group, inner and outer automorphisms
are regarded as sitting in a short exact sequence
$$
  \xymatrix{
    \mathrm{Inn}(G)
    \ar[r]
    &
    \mathrm{Aut}(G)
    \ar[r]
    &
    \mathrm{Out}(G)
  }
$$
of ordinary groups.

But we will find shortly that we ought to be regarding the 
conjugation automorphisms by two group elements which differ
by an element in the center of the group as \emph{different}
inner automorphisms. 

So adopting this point of view for ordinary groups, 
one gets instead the exact sequence
$$
  \xymatrix{
    \mathrm{Z}(G)
    \ar[r]
    &    
    \mathrm{Inn}'(G)
    \ar[r]
    &
    \mathrm{Aut}(G)
    \ar[r]
    &
    \mathrm{Out}(G)
  }
  \,.
$$

Of course this means setting $\mathrm{Inn}'(G) \simeq G$, which 
seems to make this step rather ill motivated. But it turns out
that this degeneracy of concepts is a coincidence of low dimensions
and will be lifted as we pass to inner automorphisms of 
higher groups.

First recall the standard definitions of center and automorphism 
of 2-groupoids:

\begin{definition}
  Given any strict $2$-groupoid $C$,
  \begin{itemize}
    \item
      the automorphism $3$-group
      $$
        \mathrm{AUT}(C) := \mathrm{Aut}_{2\mathrm{Cat}}(C)
      $$
      is the $2$-groupoid of isomorphisms on $C$: 
      objects are the strict and strictly invertible 2-functors 
      $$
        \xymatrix{
           C 
             \ar[rr]^f_\simeq
             &&
           C
        }
        \,,
      $$ 
      morphisms are pseudonatural 
      transformations 
      $$
        \xymatrix{
          C
          \ar@/^2pc/[rr]^f_{\ }="s"
          \ar@/_2pc/[rr]_{f'}^{\ }="t"
          &&
          C
          \ar@{=>}^h "s"; "t"
        }
      $$
      between these and 2-morphisms are modifications
\[
\xymatrix{
  \makebox(8,8){}
  \ar@/^2pc/[rr]|{h}
  \ar@{-->}@/_2pc/[rr]|{\tiny \mbox{$h'$}}
  &&
  \makebox(8,8){}
  \ar@{==>}^\rho (12,3); (12,-3)
  \ar@<+8pt>@/^2.4pc/|<<<<<<<<<<<<{f'} (12,11.5); (12,-14.5)
  \ar@<-8pt>@/_2.4pc/|<<<<<<<<<<<<{f} (12,11.5); (12,-14.5)
  \ar@{}|{\small \mbox{$C$}} (13,-17.0); (13,-17.0)
  \ar@{}|{\tiny \mbox{$C$}} (12,12); (12,12)
}
\]
between those --
the product on the 3-group comes from the composition of autofunctors;
    \item
      the center of $C$
      $$
        Z(C) := \mathbf{B} \mathrm{AUT}(\mathrm{Id}_C)
      $$
      is the (suspended) automorphism 2-group of the identity
      on $C$, i.e. the full subgroupoid of $\mathrm{AUT}(C)$ on the single
      object $\mathrm{Id}_C$.
  \end{itemize}
\end{definition}

\paragraph{Example.}
The automorphism 2-group of any ordinary group $G$
(regarded as a 2-group $\mathrm{Disc}(G)$ with only identity morphisms)
$$
  \mathrm{AUT}(G) := \mathrm{AUT}(\mathbf{B} G)
$$
is that coming from the crossed module
$$
  \xymatrix{
    G \ar[r]^<<<<<{\mathrm{Ad}} & \mathrm{Aut}(G)
    \ar[r]^{\mathrm{Id}}
    &
    \mathrm{Aut}(G)
    \,.
  }
$$
The center 
$$
  Z(G) := Z(\mathbf{B} G)
	$$
of any ordinary group is indeed the ordinary center of the
group, regarded as a 1-object category.

\paragraph{Example.} The automorphism 3-group of a strict 2-group (conceived in terms
of crossed modules and 2-crossed modules) is discussed in theorem 4.3 of 
\cite{BrownIcen}.

To these two standard definitions, we add the following one, 
which is supposed to be the proper $2$-categorical generalization
of the concept of inner automorphisms.

\begin{definition}[inner automorphisms]
  Given any strict $2$-groupoid $C$, the tangent $2$-groupoid
  $$
    \mathrm{INN}(C) := T_{\mathrm{Id}_C}(\mathrm{Aut}_{2\mathrm{Cat}}(C))
  $$ 
  is called the $2$-groupoid of inner automorphisms of $C$, and as 
  such thought of as being equipped with the monoidal structure
  inherited from $\mathrm{End}(C)$.
\end{definition}
If the transformation starting at the identity is denoted $q$,
it makes good sense to call the inner automorphism being the
target of that transformation $\mathrm{Ad}_q$:
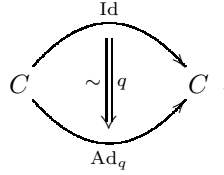
$$
    \xymatrix{
      C
      \ar@/^1.9pc/[rr]^{\mathrm{Id}}_{\ }="s"
      \ar@/_1.9pc/[rr]_{\mathrm{Ad}_q}^{\ }="t"
      &&
      C
      \ar@{=>}^q_\sim "s"; "t"
    }
    \,.
$$
A bigon of this form is an object in $\mathrm{INN}(C)$. The product of
two such objects is the horizontal composition of these bigons
in $2\mathrm{Cat}$. We shall spell this out in great detail for the
case $C = \mathbf{B} G_{(2)}$ in \ref{The 3-group inng2}.

\begin{proposition}
  For $C$ any strict 2-category, 
  we have canonical morphisms
  $$
    \xymatrix{
      \mathrm{Z}(C)
      \ar@{^{(}->}[r]
      &
      \mathrm{INN}(C)
      \ar[r]
      &
      \mathrm{AUT}(C)
    }
  $$
  of strict 2-categories
  whose composition sends everything to the identity 2-morphism
  on the identity 1-morphism on the identity automorphism of $C$.
  
Moreover, this 
sits inside the exact sequence from proposition \ref{exact sequence}
as 
$$
  \raisebox{20pt}{
  \xymatrix{
     Z(C)
     \ar[r]
     \ar@{^{(}->}[d]
     &
     \mathrm{INN}(C)
     \ar[r]
     \ar@{^{(}->}[d]
     &
     \mathrm{AUT}(C)
     \ar@{^{(}->}[d]
     \\
     \mathrm{Mor}(\mathbf{C})
     \ar[r]
     &
     T \mathbf{C}
     \ar[r]
     &
     \mathbf{C}
  }
  }
  \,,
$$
where $\mathbf{C} := \mathrm{Aut}_{2\mathrm{Cat}}(C)$.
\end{proposition}
\proof
  Recall that a morphism in $Z(C)$ is a transformation
  of the form
$$
    
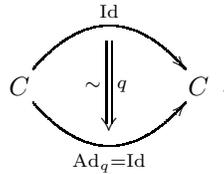
\xymatrix{
      C
      \ar@/^1.9pc/[rr]^{\mathrm{Id}}_{\ }="s"
      \ar@/_1.9pc/[rr]_{\mathrm{Ad}_q = \mathrm {Id}}^{\ }="t"
      &&
      C
      \ar@{=>}_\sim^q "s"; "t"
    }
    \,.
$$
This gives the obvious inclusion $Z(G) \hookrightarrow \mathrm{INN}(G)$.
The morphism $\mathrm{INN}(G) \to \mathrm{AUT}(G)$
maps
  $$
    \xymatrix{
      C
      \ar@/^1.9pc/[rr]^{\mathrm{Id}}_{\ }="s"
      \ar@/_1.9pc/[rr]_{\mathrm{Ad}_q}^{\ }="t"
      &&
      C
      \ar@{=>}^q "s"; "t"
    }
    \hspace{5pt}
      \mapsto
    \hspace{5pt}
   \xymatrix{
     C
     \ar[rr]^{\mathrm{Ad}_q}
     &&
     C
   }
   \,.
  $$
\endofproof

\paragraph{Remark.} One would now want to define and construct the
cokernel $\mathrm{OUT}(C)$ of the morphism 
$\mathrm{INN}(C) \to \mathrm{AUT}(C)$ and then say that
\begin{eqnarray}
 \label{exact sequence including outers}
 \xymatrix{
     Z(C)
     \ar[r]
     &
     \mathrm{INN}(C)
     \ar[r]
     &
     \mathrm{AUT}(C)  
     \ar[r]
     &
     \mathrm{OUT}(C)
  }
\end{eqnarray}
is an exact sequence of 3-groups. 
We shall not consider $\mathrm{OUT}(C)$ here for proper 3-groups.
Restricted to just 2-groups, however, one obtains the situation
described in the next section.

\subsection{Inner automorphism $2$-groups}

Even though inner automorphism 2-groups of ordinary (1-)groups are just
a special case of the inner automorphism 3-groups of strict 2-groups to be
discussed in the following, it may be helpful to spell out this simpler
case in detail, in order to see how it connects with familiar
examples of crossed modules.

For $G$ any (ordinary) group, the sequence \ref{exact sequence including outers}
is the exact sequence of 2-groups
$$
  \hspace{-.2cm}
  \xymatrix@R=11pt@C=12pt{
    1 
    \ar[r]
    &
    Z(G)
    \ar[r]
    &
    \mathrm{INN}(G)
    \ar[r]
    &
    \mathrm{AUT}(G)
    \ar[r]
    &
    \mathrm{OUT}(G)
    \ar[r]
    &
    1
    \\
    (1 \to 1) 
    \ar[r]
    &
    \ar[r]
    (1 \to Z(G))
    \ar[r]
    &
    (G \stackrel{\mathrm{Id}}{\to} G)
    \ar[r]
    &
    (G \stackrel{\mathrm{Ad}}{\to} \mathrm{Aut}(G))
    \ar[r]
    &
    (1 \to  \mathrm{Out}(G))    
    \ar[r]
    &
    (1 \to 1)
  }
$$
corresponding to the exact sequence of crossed modules given in the second line.

Notice that there exists also the crossed module $(\mathrm{Inn}(G) \to \mathrm{Aut}(G))$,
which however does not appear in the above sequence. In particular, this crossed module is not
the one corresponding to our 2-group $\mathrm{INN}(G)$, which we had described in detail
in the introduction.

\subsection{Inner automorphism $3$-groups.}

Now we apply the general concept of inner automorphisms
to $2$-groups. The following definition just establishes
the appropriate shorthand notation.
\begin{definition}
  For $G_{(2)}$ a strict $2$-group, we write
  $$
    \mathrm{INN}(G_{(2)}) := \mathrm{INN}(\mathbf{B} G_{(2)})
  $$
  for its $3$-group of inner automorphisms.
\end{definition}
In general this notation could be ambiguous, since one might
want to consider the inner automorphisms of just the
1-groupoid underlying $G_{(2)}$. However, in the present
context this will never occur and using the above definition
makes a couple of expressions more manifestly appear as generalizations
of familiar ones.

\paragraph{Example.}
For $G$ an ordinary group, regarded as a discrete 2-group, one finds that
$$
  \mathrm{INN}(G) := 
    T_{\mathrm{Id}_{\mathbf{B} G}}(n\mathrm{Cat}) 
    \simeq T_\bullet (\mathbf{B} G)
$$
is the codiscrete groupoid over the elements of $G$. 
Its nature as a groupoid is manifest from its realization
as
$$
  \mathrm{INN}(G) = T_\bullet (\mathbf{B} G)
  \,.
$$
But it is also a (strict) 2-group. The monoidal structure is that
coming from its realization as
$
 \mathrm{INN}(G) := T_{\mathrm{Id}_{\mathbf{B} G}}(n\mathrm{Cat})
 \,.
$
The crossed module corresponding to this strict 2-group is
$$
  \xymatrix{
    G \ar[r]^{\mathrm{Id}} & G \ar[r]^<<<<<<{\mathrm{Ad}} & \mathrm{Aut}(G)
    \,.
  }
$$
The main point of interest for us is the generalization of this
fact to strict 2-groups. One issue that one needs to be aware of
then is that the above 
isomorphism
$    T_{\mathrm{Id}_{\mathbf{B} G}}(n\mathrm{Cat}) 
    \simeq T_\bullet (\mathbf{B} G)$ becomes a mere inclusion.

\begin{proposition}
\label{the embedding}
For $G_{(2)}$ any strict 2-group, we have an inclusion
$$
  T_{\bullet} \mathbf{B} G_{(2)}
  \subset
  T_{\mathrm{Id}_{\mathbf{B} G_{(2)}}}(\mathrm{Aut}_{2\mathrm{Cat}}(\mathbf{B} G_{(2)}))
$$
of strict 2-groupoids.

This realizes $T_{\bullet} \mathbf{B} G_{(2)}$ as a sub 2-groupoid of
$T_{\mathrm{Id}_{\mathbf{B} G_{(2)}}}(\mathrm{Aut}_{2\mathrm{Cat}}(\mathbf{B} G_{(2)}))$.
\end{proposition}
\proof
    The inclusion is essentially fixed by its action on objects:
    we define that an object in $T_\bullet \mathbf{B} G_{(2)}$, which is a morphism
    $$
      \xymatrix{
        \bullet
        \ar[r]^q
        &
        \bullet
      }
    $$
    in $\mathbf{B} G$, is sent to the conjugation automorphism
\[
\begin{array}{ccccc}
\mathrm{Ad}_{q}
\hspace{10pt}
&:&
\hspace{10pt}
\mathbf{B} G_{(2)}
&\rightarrow&
\mathbf{B} G_{(2)}
\\
\\
\hspace{10pt}
&&
\hspace{10pt}
\xymatrix{
  \bullet
    \ar@/^1.8pc/[rr]^{g}_{\ }="s"
    \ar@/_1.8pc/[rr]_{g'}^{\ }="t"
    &&
  \bullet
  \ar@{=>}^{h} "s"; "t"
}
\hspace{10pt}
 & \mapsto &
\hspace{10pt}
\xymatrix{
  \bullet
  \ar[r]^{q^{-1}}
  &
  \bullet
    \ar@/^1.8pc/[rr]^{g}_{\ }="s"
    \ar@/_1.8pc/[rr]_{g'}^{\ }="t"
    &&
  \bullet
  \ar[r]^{q}
  &
  \bullet
  \ar@{=>}^{h} "s"; "t"
}
\end{array}
\]
The transformation
$$
    \xymatrix{
      \mathbf{B} G_{(2)}
      \ar@/^1.9pc/[rr]^{\mathrm{Id}}_{\ }="s"
      \ar@/_1.9pc/[rr]_{\mathrm{Ad}_q}^{\ }="t"
      &&
      \mathbf{B} G_{(2)}
      \ar@{=>}^q_\sim "s"; "t"
    }
    \,.
$$
connecting this to the identity is given by the component map
$$
  (
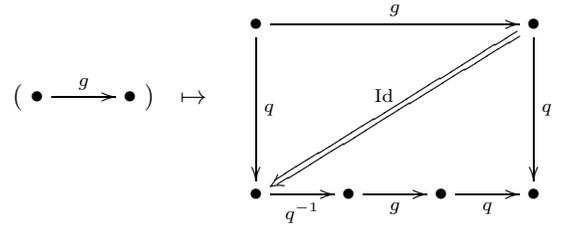
\xymatrix{
    \bullet 
    \ar[r]^g
    &
    \bullet
  })
  \hspace{7pt}
    \mapsto
  \hspace{7pt}
  \raisebox{30pt}{
  \xymatrix{
    \bullet
    \ar[dd]^q
    \ar[rrr]^g
    &&&
    \bullet
    \ar[dd]^{q}
    \ar@{=>}[ddlll]_{\mathrm{Id}}
    \\
    \\
    \bullet 
    \ar[r]_{q^{-1}}
    &
    \bullet
    \ar[r]_{g}    
    &
    \bullet
    \ar[r]_{q}    
    &
    \bullet
  }
  }
  \,.
$$
In general one could consider transformations whose component
maps involve here a non-identity 2-morphism. The inclusion
we are describing picks out excactly those transformations
whose component map only involves identity 2-morphisms.

The crucial point to realize now is the
form of the component maps of morphisms 
$$
  
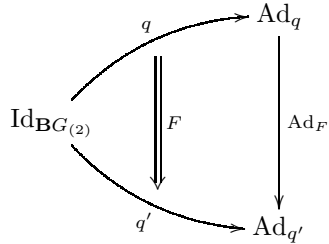
\xymatrix{
   && \mathrm{Ad}_q
   \ar[dd]^{\mathrm{Ad}_F}
   \\
   \mathrm{Id}_{\mathbf{B} G_{(2)}}
   \ar@/^1pc/[urr]^{q}_{\ }="s"
   \ar@/_1pc/[drr]_{q'}^{\ }="t"
   \\
   && \mathrm{Ad}_{q'}
   \ar@{=>}^F "s"; "t"
  }
$$
in 
$T_{\mathrm{Id}_{\mathbf{B} G_{(2)}}}(\mathrm{Aut}_{2\mathrm{Cat}}(\mathbf{B} G_{(2)}))$.

The corresponding component map equation is
$$
  \raisebox{60pt}
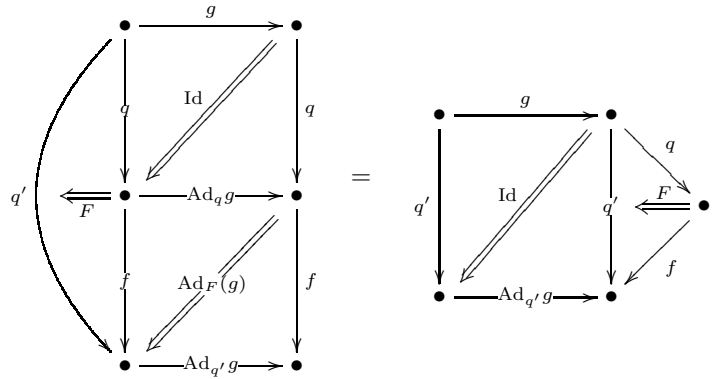
{
  \xymatrix{
    \bullet
    \ar[rr]^g_>{\ }="s1"
    \ar[dd]|{q}
    \ar@/_2.8pc/[dddd]_{q'}^{\ }="t3"
    &&
    \bullet
    \ar[dd]^q
    \\
    \\
    \bullet
    \ar[rr]|{\mathrm{Ad}_q g}^<{\ }="t1"
    \ar[dd]|{f}
    \ar@{}|{\ }="s3"
    &&
    \bullet
    \ar[dd]^{f}_<{\ }="s2"
    \\
    \\
    \bullet
    \ar[rr]|{\mathrm{Ad}_{q'} g}^<{\ }="t2"
    &&
    \bullet    
    \ar@{=>}_{\mathrm{Id}} "s1"; "t1"
    \ar@{=>}|{\mathrm{Ad}_F(g)} "s2"; "t2"
    \ar@{=>}^F "s3"; "t3"
  }
  }
  \hspace{7pt}
    =
  \hspace{7pt}
  \raisebox{26pt}{
  \xymatrix{
    \bullet
    \ar[rr]^g_>{\ }="s1"
    \ar[dd]_{q'}
    &&
    \bullet
    \ar[dd]|{q'}^{\ }="t2"
    \ar[dr]^{q}
    \\
    &&& 
    \bullet
    \ar@{}|{\ }="s2"
    \ar[dl]^{f}
    \\
    \bullet
    \ar[rr]|{\mathrm{Ad}_{q'} g}^<{\ }="t1"
    &&
    \bullet
    \ar@{=>}_{\mathrm{Id}} "s1"; "t1"
    \ar@{=>}_{F} "s2"; "t2"
  }
  }
  \,.
$$
Solving this for $\mathrm{Ad}_F$ shows that this 
is given by conjugation
$$
 \mathrm{Ad}_F
  \hspace{8pt}
   :
  \hspace{8pt}
  (\xymatrix{
    \bullet
    \ar[r]^g
    &
    \bullet
  })
  \hspace{7pt}
   \mapsto
\hspace{7pt}
  \raisebox{36pt}
  {
   \xymatrix{
     \bullet
     \ar@/^1pc/[drr]^{q^{-1}}_{\ }="s_1"
     \ar[dd]_{f}
     &&
     &&
     &&
     \bullet
     \ar[dd]^{f}
  \\
  &&
     \bullet
     \ar[rr]^g
     &&
     \bullet
     \ar@/^1pc/[urr]^{q}_{\ }="s_2"
     \ar@/_1pc/[drr]_{q'}^{\ }="t_2"
     &&
  \\
     \bullet
     \ar@/_1pc/[urr]_{q'^{-1}}^{\ }="t_1"
     &&
     &&
     &&
     \bullet
     \ar@{=>} "s_1";"t_1"
     \ar@{=>}^{F} "s_2";"t_2"
  }
  }
$$
with a morphism in $T_\bullet (\mathbf{B} G_{(2)})$.
And each such morphism in 
$T_\bullet (\mathbf{B} G_{(2)})$
yields a morphism in 
$T_{\mathrm{Id}_{\mathbf{B} G_{(2)}}}(\mathrm{Aut}_{2\mathrm{Cat}}(\mathbf{B} G_{(2)}))$
this way.

Finally, 2-morphisms in 
$T_{\mathrm{Id}_{\mathbf{B} G_{(2)}}}(\mathrm{Aut}_{2\mathrm{Cat}}(\mathbf{B} G_{(2)}))$
between these morphisms
\[
\xymatrix{
  \makebox(8,8){}
  \ar@/^2pc/[rr]|{\mathrm{Ad}_F}
  \ar@{-->}@/_2pc/[rr]|{\tiny \mbox{$\mathrm{Ad}_K$}}
  &&
  \makebox(8,8){}
  \ar@{==>}^L (12,3); (12,-3)
  \ar@<+8pt>@/^2.4pc/|<<<<<<<<<<<<{\mathrm{Ad}_{q'}} (12,11.5); (12,-14.5)
  \ar@<-8pt>@/_2.4pc/|<<<<<<<<<<<<{\mathrm{Ad}_{q}} (12,11.5); (12,-14.5)
  \ar@{}|{\small \mbox{$\mathbf{B} G_{(2)}$}} (13,-17.0); (13,-17.0)
  \ar@{}|{\tiny \mbox{$\mathbf{B} G_{(2)}$}} (12,12); (12,12)
}
\]
come from component maps
\[
  \bullet 
  \hspace{6pt}
  \mapsto
  \hspace{10pt}
  \xymatrix{
    \bullet
    \ar@/^1.8pc/[rr]^{f}
    \ar@/_1.8pc/[rr]_{k}
    &&
    \bullet
    \ar@{=>}^{L} (11,3); (11,-3)
  }
  \hspace{3pt}
  \in \hspace{3pt}
  \Mor_2\of{\mathbf{B} G_{(2)}}
  \,.
\]
A sufficient condition for these component maps to solve the
required condition for modifications of pseudonatural transformations is
that they make
$$
\raisebox{36pt}{
\xymatrix{
&&\bullet \ar @/^1pc/ [dd]^{k}_{}="t_2"\ar @/_1pc/ [dd]_{f}^{}="s_2"
\\
\bullet \ar @/^1pc/ [urr]^{q}_{\ }="s" \ar @/_1pc/ [drr]^{\ }="t"_{}&&
\\
&&\bullet 
\ar @{=>} @/_.7pc/ "s";"t"_{F} |>\hole
\ar @{:>} @/^.7pc/ "s";"t"_{K} |>\hole
\ar @2 "s_2";"t_2"^{L}
}
}
$$
2-commute. 
But this defines a 2-morphism in $T_\bullet \mathbf{B} G_{(2)}$.
And each such 2-morphism in 
$T_\bullet (\mathbf{B} G_{(2)})$
yields a 2-morphism in 
$T_{\mathrm{Id}_{\mathbf{B} G_{(2)}}}(\mathrm{Aut}_{2\mathrm{Cat}}(C))$
this way.
\endofproof

The crucial point is that by the embedding
$$
  T_\bullet \mathbf{B} G_{(2)}
  \subset
  T_{\mathrm{Id}_{\mathbf{B} G_{(2)}}}(\mathrm{Aut}_{2\mathrm{Cat}}(\mathbf{B} G_{(2)}))
$$
the former 2-category inherits the monoidal structure of the latter
and hence becomes a 3-group in its own right. This 3-group is
the object of interest here.

\section{The 3-group $\mathrm{INN}_0(G_{(2)})$}
\label{The 3-group inng2}

\begin{definition}[$\mathrm{INN}_0(G_{(2)})$]
  For $G_{(2)}$ any strict 2-group, 
  the 3-group $\mathrm{INN}_0(G_{(2)})$ is, as a 2-groupoid,
  given by
  $$
    \mathrm{INN}_0(G_{(2)}) := T_\bullet \mathbf{B} G_{(2)}
  $$
  and equipped with the monoidal structure inherited from the
  embedding of proposition \ref{the embedding}.
\end{definition}

We now describe $\mathrm{INN}_0(G_{(2)})$ for $G_{(2)}$
coming from the crossed module
  $$
  \xymatrix{
    H \ar[r]^t & G \ar[r]^(.35)\alpha & \mathrm{Aut}(H)
  }
  $$
in more detail, in particular spelling 
out the monoidal structure.
We extract the operations in the
crossed module corresponding to the various
compositions in $\mathrm{INN}_0(G_{(2)})$ and then
finally identify the 2-crossed module encoded by this.

\subsection{Objects}

The objects of $\mathrm{INN}_0(G_{(2)})$ are exactly the objects of $G_{(2)}$,
hence the elements of $G$:
$$
  \mathrm{Obj}(\mathrm{INN}(G_{(2)}))
  =
  G
  \,.
$$

The product of two objects in 
$\mathrm{INN}(G_{(2)})$ is just the product in $G$.

\subsection{Morphisms}

The morphisms 
$$
  g \to h
$$
in $\mathrm{INN}(G_{(2)})$  are
\begin{eqnarray*}
\mathrm{Mor}(\mathrm{INN}_0(G_{(2)}))
&=&
\left\lbrace
\left.
\raisebox{35pt}{
\xymatrix{
&&\bullet\ar[dd]^{f}\\
\bullet \ar @/^1pc/ [urr]^{g}_{\ }="s" \ar @/_1pc/ [drr]^{\ }="t"_{h}&&\\
&&\bullet 
\ar @2 "s";"t"^{F}
}
}
\hspace{7pt}
\right|
\hspace{7pt}
{
{f,g,h \in G, \; F \in H}
\atop
{
  h = t(F)fg
}
}
\right\rbrace
\\
\\
\\
&=&
\lbrace
  (f,F;g)
  \;
  |
  \;
  f,g \in G, F \in H
\rbrace
\,.
\end{eqnarray*}

\subsubsection{Composition}

The composition of two such morphisms
$$
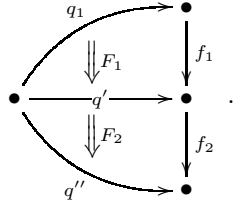
\raisebox{36pt}{
\xymatrix{
&&
 \bullet\ar[d]^{f_1}
 \\
 \bullet 
  \ar @/^1pc/ [urr]^{q_1}_{\ }="s_1" 
  \ar @/_1pc/ [drr]^{\ }="t_2"_{q''} 
  \ar[rr]^(.45){\ }="t_1"_(.45){\ }="s_2" |{q'}&&\bullet \ar[d]^{f_2}\\
 && 
 \bullet
 \ar @2 "s_1";"t_1"^{F_1}
 \ar @2 "s_2";"t_2"^{F_2}
}}
\,.
$$
%%%%%%%%%%%%%%%%%%%%%%
is in terms of group labels given by
$$
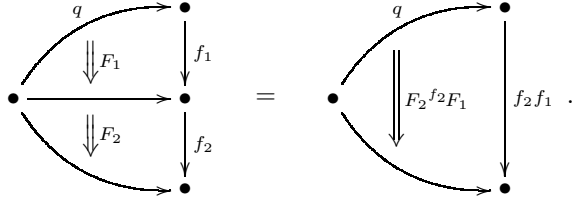
\raisebox{36pt}
{
\xymatrix{
&&\bullet \ar[d]^{f_1}\\
 \bullet 
  \ar @/^1pc/ [urr]^{q}_{\ }="s_1" 
  \ar @/_1pc/ [drr]^{\ }="t_2"_{} 
  \ar[rr]^(.45){\ }="t_1"_(.45){\ }="s_2" 
%|{q_2t(F_1)f_1q_1}     %Label for middle arrow - to big to fit in break
&&\bullet\ar[d]^{f_2}\\
&&\bullet 
\ar @2 "s_1";"t_1"^{F_1}
\ar @2 "s_2";"t_2"^{F_2}
}
}
\equals
\raisebox{36pt}
{
\xymatrix{
&&\bullet\ar[dd]^{\gprod{f_1}{f_2}}\\
\bullet 
 \ar @/^1pc/ [urr]^{q}_(.4){\ }="s" 
 \ar @/_1pc/ [drr]^(.4){\ }="t"_{}&&\\
&&\bullet 
\ar @2 "s";"t"^{\hprod{\rightact{f_2}{F_1}}{F_2}}
}
}
\,.
$$

\subsubsection{Product}

Horizontal composition of automorphisms $\mathbf{B} G_{(2)} \to \mathbf{B} G_{(2)}$
gives the product in the 3-group $\mathrm{INN}(G_{(2)})$

Left whiskering of pseudonatural transformations
$$
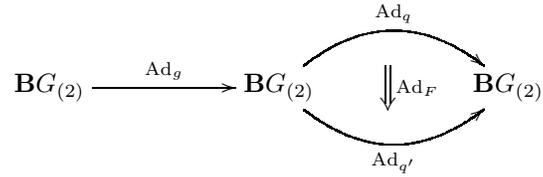
\xymatrix{
  \mathbf{B} G_{(2)}
    \ar[rr]^{\mathrm{Ad}_{g}}
    &&
  \mathbf{B} G_{(2)}
    \ar@/^1.8pc/[rr]^{\mathrm{Ad}_{q}}
    \ar@/_1.8pc/[rr]_{\mathrm{Ad}_{q'}}
    &&
  \mathbf{B} G_{(2)}
  \ar@{=>}^{\mathrm{Ad}_F} (19,3)+(26,0); (19,-3)+(26,0)
}
$$
amounts to the operation
$$
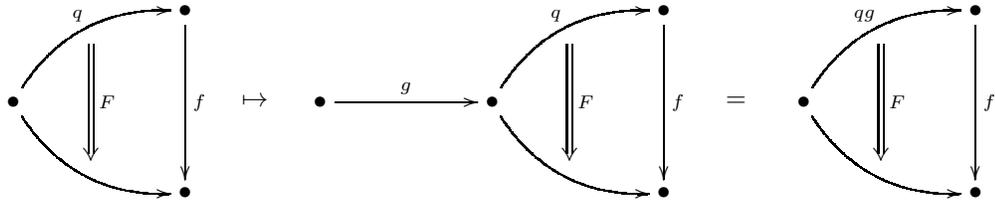
\raisebox{36pt}
{
\xymatrix{
&&
 \bullet\ar[dd]^{f}
  \\
  \bullet 
  \ar @/^1pc/ [urr]^{q}_{\ }="s" \ar @/_1pc/ [drr]^{\ }="t"_{}&&\\
&&\bullet 
\ar @2 "s";"t"^{F}
}
}
\hspace{8pt}
  \mapsto
\hspace{8pt}
\raisebox{36pt}
{
\xymatrix{
&&&&\bullet\ar[dd]^{f}\\
\bullet\ar[rr]^g&&\bullet \ar @/^1pc/ [urr]^{q}_{\ }="s" \ar @/_1pc/ [drr]^{\ }="t"_{}&&\\
&&&&\bullet 
\ar @2 "s";"t"^{F}
}
}
\equals
\raisebox{36pt}
{
\xymatrix{
&&\bullet\ar[dd]^{f}\\
\bullet \ar @/^1pc/ [urr]^{\gprod{g}{q}}_{\ }="s" 
\ar @/_1pc/ [drr]^{\ }="t"_{}&&\\
&&\bullet 
\ar @2 "s";"t"^{\leftact{g}{F}}
}
}
$$
on the corresponding triangles.

Right whiskering of pseudonatural transformations
$$
\xymatrix{
  \mathbf{B} G_{(2)}
    \ar@/^1.8pc/[rr]^{\mathrm{Ad}_{q}}
    \ar@/_1.8pc/[rr]_{\mathrm{Ad}_{q'}}
    &&
  \mathbf{B} G_{(2)}
    \ar[rr]^{\mathrm{Ad}_{q}}
    &&
  \mathbf{B} G_{(2)}
  \ar@{=>}^{\mathrm{Ad}_F} (14,3); (14,-3)
}
$$
amounts to the operation
$$
\raisebox{36pt}
{
\xymatrix{
&&
 \bullet\ar[dd]^{f}
  \\
  \bullet 
  \ar @/^1pc/ [urr]^{q}_{\ }="s" \ar @/_1pc/ [drr]^{\ }="t"_{}&&\\
&&\bullet 
\ar @2 "s";"t"^{F}
}
}
\hspace{8pt}
  \mapsto
\hspace{8pt}
\raisebox{36pt}
{
\xymatrix{
&&\bullet \ar[dd]^{f} \ar[rr]^{g} \ar@{}[ddrr] |{=} && \bullet 
 \ar@{-->}[dd]^{\gprod{g^{-1}}{\gprod{f}{g}}}\\
\bullet 
 \ar @/^1pc/ [urr]^{q}_{\ }="s" 
\ar @/_1pc/ [drr]^{\ }="t"_{}&&&&\\
&&
 \bullet \ar[rr]_g &&\bullet
\ar @2 "s";"t"^{F}
}
}
\equals
\raisebox{36pt}
{
\xymatrix{
&&\bullet\ar[dd]^{\gprod{g^{-1}}{\gprod{f}{g}}}\\
\bullet 
  \ar @/^1pc/ [urr]^{\gprod{q}{g}}_{\ }="s" 
 \ar @/_1pc/ [drr]^{\ }="t"_{}&&\\
&&\bullet 
\ar @2 "s";"t"^{\rightact{g}{F}}
}
}
$$
on the corresponding triangles.

Since $2\mathrm{Cat}$ is a Gray-category, 
the horizontal composition of pseudonatural transformations
\[
\xymatrix{
  \mathbf{B} G_{(2)}
    \ar@/^1.8pc/[rr]^{\mathrm{Ad}_{q_1}}
    \ar@/_1.8pc/[rr]_{\mathrm{Ad}_{q'_1}}
    &&
  \mathbf{B} G_{(2)}
    \ar@/^1.8pc/[rr]^{\mathrm{Ad}_{q_2}}
    \ar@/_1.8pc/[rr]_{\mathrm{Ad}_{q'_2}}
    &&
  \mathbf{B} G_{(2)}
  \ar@{=>}^{\mathrm{Ad}_{F_1}} (14,3); (14,-3)
  \ar@{=>}^{\mathrm{Ad}_{F_2}} (19,3)+(26,0); (19,-3)+(26,0)
}
\]
is ambiguous. We shall agree to read this as
$$
\raisebox{20pt}{
\xymatrix{
  \mathbf{B} G_{(2)}
    \ar@/^1.8pc/[rr]^{\mathrm{Ad}_{q_1}}
    \ar@/_1.8pc/[rr]_{\mathrm{Ad}_{q'_1}}
    &&
  \mathbf{B} G_{(2)}
    \ar@/^1.8pc/[rr]^{\mathrm{Ad}_{q_2}}
    &&
  \mathbf{B} G_{(2)}
  \ar@{=>}^{\mathrm{Ad}_{F_1}} (14,3); (14,-3)
  \\
  \mathbf{B} G_{(2)}
    \ar@/_1.8pc/[rr]_{\mathrm{Ad}_{q'_1}}
    &&
  \mathbf{B} G_{(2)}
    \ar@/^1.8pc/[rr]^{\mathrm{Ad}_{q_2}}
    \ar@/_1.8pc/[rr]_{\mathrm{Ad}_{q'_2}}
    &&
  \mathbf{B} G_{(2)}
  \ar@{=>}^{\mathrm{Ad}_{F_2}} (19,3)+(26,-15); (19,-3)+(26,-15)
}
}
\,.
$$
The corresponding operation on triangles labelled in the crossed
module is
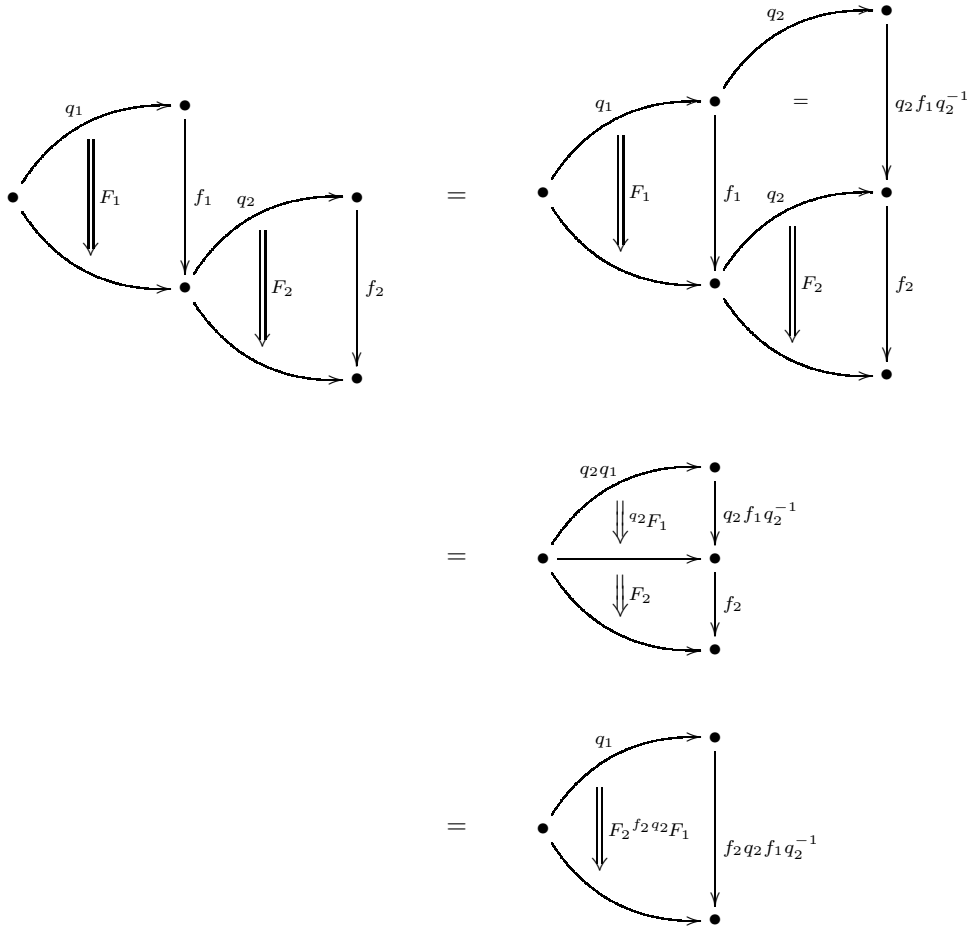
\begin{eqnarray*}
\raisebox{36pt}
{
\xymatrix{
&&\bullet \ar[dd]^{f_1}                                                                          &&\\
\bullet \ar @/^1pc/ [urr]^{q_1}_{\ }="s_1" \ar @/_1pc/ [drr]^{\ }="t_1"_{}&& &&\bullet\ar[dd]^{f_2}\\
&&\bullet \ar @/^1pc/ [urr]^{q_2}_{\ }="s_2" \ar @/_1pc/ [drr]^{\ }="t_2"_{}
\ar @2 "s_1";"t_1"^{F_1}
\\
&&&&\bullet
\ar @2 "s_2";"t_2"^{F_2}
}
}
&\equals&
\raisebox{72pt}
{
\xymatrix{
&&&&
\bullet
\ar[dd]^{
  \gprod{q_2^{-1}}{\gprod{f_1}{q_2}}
}\\
&&\bullet \ar[dd]^{f_1}   \ar @/^1pc/ [urr]^{q_2}  \ar@{}[rr] |{=}                                                     &&\\
\bullet  
\ar @/^1pc/ [urr]^{q_1}_{\ }="s_1" 
\ar @/_1pc/ [drr]^{\ }="t_1"_{}&& &&\bullet\ar[dd]^{f_2}\\
&&
\bullet 
\ar @/^1pc/ [urr]^{q_2}_{\ }="s_2" 
\ar @/_1pc/ [drr]^{\ }="t_2"_{}
\ar @2 "s_1";"t_1"^{F_1}
\\
&&&&\bullet
\ar @2 "s_2";"t_2"^{F_2}
}
}
\\ %End first line
&&\\ %Gap
&&\\ %Gap
&\equals& 
\raisebox{36pt}
{
\xymatrix{
&&
\bullet 
\ar[d]^{ 
  \gprod{q_2^{-1}}{\gprod{f_1}{q_2}}
}
\\
\bullet \ar @/^1pc/ [urr]^{\gprod{q_1}{q_2}}_{\ }="s_1"
\ar @/_1pc/ [drr]^{\ }="t_2"_{} 
\ar[rr]^(.45){\ }="t_1"_(.45){\ }="s_2" %|{q_2t(F_1)f_1q_1}
&&\bullet\ar[d]^{f_2}\\
&&\bullet 
\ar @2 "s_1";"t_1"^{\rightact{q_2}{F_1}}
\ar @2 "s_2";"t_2"^{\leftact{q'_1}{F_2}}
}
}
\\ %End second line
&&\\ %Gap
&&\\ %Gap
&\equals& 
\raisebox{36pt}
{
\xymatrix{
&&
\bullet
\ar[dd]^(.6){
 \gprod{\gprod{q_2^{-1}}{\gprod{f_1}{q_2}}}{f_2}
}\\
\bullet \ar @/^1pc/ [urr]^{q_1}_(.35){\ }="s" 
\ar @/_1pc/ [drr]^(.35){\ }="t"_{}&&\\
&&\bullet 
\ar @2 "s";"t"^{
    \hprod{\rightact{\gprod{q_2}{f_2}}{F_1}}{F_2}
               }
}
}
\end{eqnarray*}

The non-identitcal isomorphism which relates this to the other
possible way to evaluate the horizontal composition of
pseudonatural transformations gives rise to the Peiffer
lifting of the corresponding 2-crossed module. This is
discussed in \ref{The corresponding 2-crossed module}.

\subsection{2-Morphisms}

The 2-morphisms in $\mathrm{INN}_0(G_{(2)})$ are given by
diagrams
$$
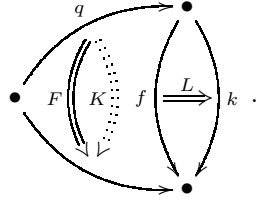
\raisebox{36pt}{
\xymatrix{
&&\bullet \ar @/^1pc/ [dd]^{k}_{}="t_2"\ar @/_1pc/ [dd]_{f}^{}="s_2"
\\
\bullet \ar @/^1pc/ [urr]^{q}_{\ }="s" \ar @/_1pc/ [drr]^{\ }="t"_{}&&
\\
&&\bullet 
\ar @{=>} @/_.7pc/ "s";"t"_{F} |>\hole
\ar @{:>} @/^.7pc/ "s";"t"_{K} |>\hole
\ar @2 "s_2";"t_2"^{L}
}
}
\,.
$$

In terms of the group labels this means that $L \in H$
satisfies
\begin{eqnarray}
  \label{relation of Inn(G_2) 2-morphisms to source and target}
  L = \leftact{q^{-1}}{\hprod{F}{K^{-1}}}
  \,.
\end{eqnarray}

\subsubsection{Composition}
\label{composition of 2-morphisms in INNG}

The horizontal composition of 2-morphisms in $\mathrm{INN}_0(G_{(2)})$
is given by 
$$
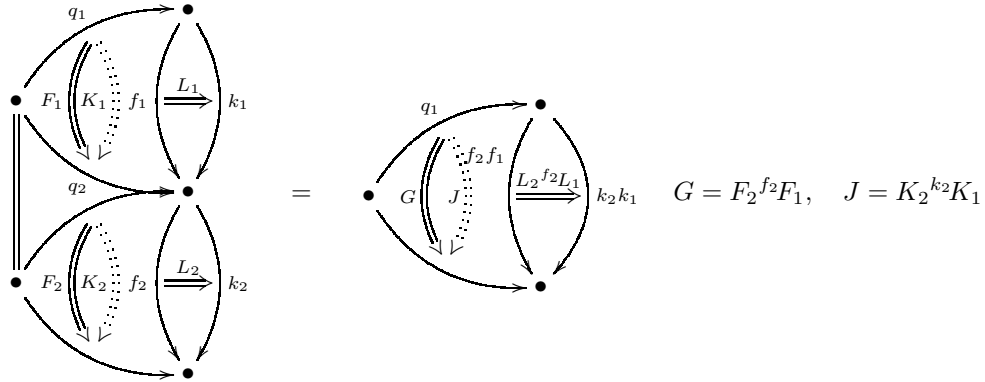
\raisebox{72pt}
{
\xymatrix{
&&\bullet \ar @/^1pc/ [dd]^{k_1}_{}="t_2"\ar @/_1pc/ [dd]_{f_1}^{}="s_2"\\
\bullet \ar @/^1pc/ [urr]^{q_1}_{\ }="s" \ar @/_1pc/ [drr]^{\ }="t"_{q_2} \ar@2{-}[dd]&&\\
&&\bullet 
\ar @{=>} @/_.7pc/ "s";"t"_{F_1} |>\hole
\ar @{:>} @/^.7pc/ "s";"t"_{K_1} |>\hole
\ar @2 "s_2";"t_2"^{L_1}
                    \ar @/^1pc/ [dd]^{k_2}_{}="t_2"\ar @/_1pc/ [dd]_{f_2}^{}="s_2"\\
\bullet \ar @/^1pc/ [urr]^{}_{\ }="s" \ar @/_1pc/ [drr]^{\ }="t"_{}&&\\
&&\bullet 
\ar @{=>} @/_.7pc/ "s";"t"_{F_2} |>\hole
\ar @{:>} @/^.7pc/ "s";"t"_{K_2} |>\hole
\ar @2 "s_2";"t_2"^{L_2}
}
}
\equals
\raisebox{36pt}
{
\xymatrix{
&&\bullet \ar @/^1.5pc/ [dd]^{\gprod{k_1}{k_2}}_{}="t_2"\ar @/_1pc/ [dd]_(.3){\gprod{f_1}{f_2}}^{}="s_2"\\
\bullet \ar @/^1pc/ [urr]^{q_1}_{\ }="s" \ar @/_1pc/ [drr]^{\ }="t"_{}&&\\
&&\bullet 
\ar @{=>} @/_.7pc/ "s";"t"_{G} |>\hole
\ar @{:>} @/^.7pc/ "s";"t"_{J} |>\hole
\ar @2 "s_2";"t_2"^{
                                \hprod{\rightact{f_2}{L_1}}{L_2}
                                }
}
}
\quad G = \hprod{\rightact{f_2}{F_1}}{F_2}, \quad J =\hprod{\rightact{k_2}{K_1}}{K_2}
$$
and
vertical composition by
$$
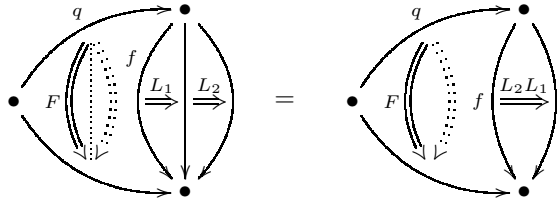
\raisebox{36pt}
{
\xymatrix{
&&\bullet \ar @/^1.5pc/ [dd]^{}_{}="t_3"  \ar[dd]^{}="s_3"_{}="t_2"       \ar @/_1.5pc/ [dd]_(.3){f}^{}="s_2"\\
\bullet \ar @/^1pc/ [urr]^{q}_{\ }="s" \ar @/_1pc/ [drr]^{\ }="t"_{}&&\\
&&\bullet 
\ar @{=>} @/_.7pc/ "s";"t"_{F} |>\hole
\ar @{:>} @/^.7pc/ "s";"t"_{} |>\hole
\ar @{.}  "s";"t"_{}
\ar @2 "s_2";"t_2"^{L_1}
\ar @2 "s_3";"t_3"^{L_2}
}
}
\equals
\raisebox{36pt}
{
\xymatrix{
&&\bullet \ar @/^1pc/ [dd]^{}_{}="t_2"\ar @/_1pc/ [dd]_{f}^{}="s_2"\\
\bullet \ar @/^1pc/ [urr]^{q}_{\ }="s" \ar @/_1pc/ [drr]^{\ }="t"_{}&&\\
&&\bullet 
\ar @{=>} @/_.7pc/ "s";"t"_{F} |>\hole
\ar @{:>} @/^.7pc/ "s";"t"_{} |>\hole
\ar @2 "s_2";"t_2"^{L_2L_1}
}
}
$$

(Notice that these compositions do go horizontally and vertically, respectively,
once we rotate such that the bigons have the standard orientation.)

Notice that whiskering along 1-morphisms
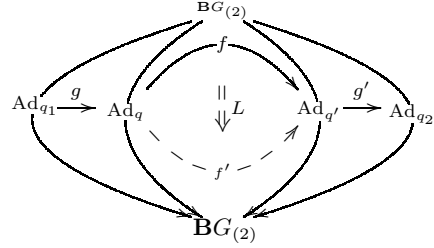
\[
\xymatrix{
  \makebox(8,8){}
  \ar@/^2pc/[rr]|{f}
  \ar@{-->}@/_2pc/[rr]|{\tiny \mbox{$f'$}}
  &&
  \makebox(8,8){}
  \ar@{==>}^L (12,3); (12,-3)
  \ar@<+8pt>@/^2.4pc/|<<<<<<<<<<<<{\mathrm{Ad}_{q'}} (12,11.5); (12,-14.5)
  \ar@<-8pt>@/_2.4pc/|<<<<<<<<<<<<{\mathrm{Ad}_{q}} (12,11.5); (12,-14.5)
  \ar@{}|{\small \mbox{$\mathbf{B} G_{(2)}$}} (12.4,-16.3); (12.4,-16.3)
  \ar@{}|{\tiny \mbox{$\mathbf{B} G_{(2)}$}} (12,13); (12,13)
  \ar@<+12pt>@/^5pc/|<<<<<<<<<<<<<<<{\mathrm{Ad}_{q_2}} (12,11.5); (12,-14.5)
  \ar@<-12pt>@/_5pc/|<<<<<<<<<<<<<<{\mathrm{Ad}_{q_1}} (12,11.5); (12,-14.5)
  \ar^{g} (-10,0); (-5,-0)
  \ar^{g'} (28,0); (33,-0)
}
\]
acts on the component maps as
$$
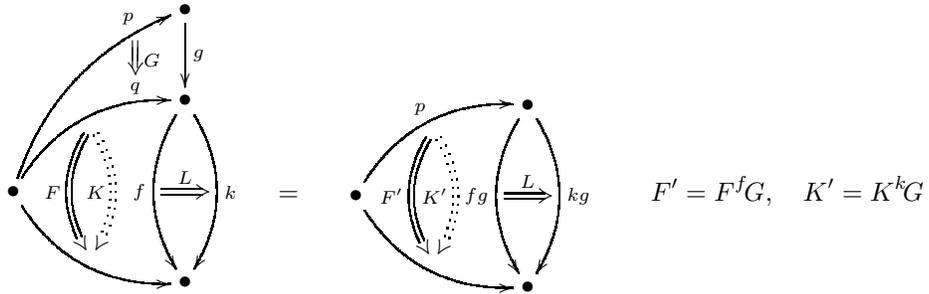
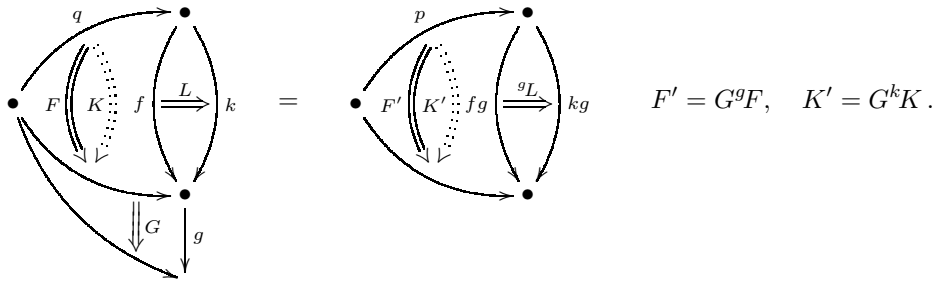
\raisebox{72pt}
{
\xymatrix{
&&\bullet \ar[d]^{g}\\
&&\bullet \ar @/^1pc/ [dd]^{k}_{}="t_2"\ar @/_1pc/ [dd]_{f}^{}="s_2"\\
\bullet \ar @/^1pc/ [urr]^(.78){q}="t_3"_{\ }="s" \ar @/_1pc/ [drr]^{\ }="t"_{}   \ar@/^1pc/ [uurr]^(.8){p}_(.8){}="s_3"&&\\
&&\bullet 
\ar @{=>} @/_.7pc/ "s";"t"_{F} |>\hole
\ar @{:>} @/^.7pc/ "s";"t"_{K} |>\hole
\ar @2 "s_2";"t_2"^{L}
\ar @{=>} "s_3";"t_3"^{G}
}
}
\equals
\raisebox{36pt}
{
\xymatrix{
&&\bullet \ar @/^1pc/ [dd]^{kg}_{}="t_2"\ar @/_1pc/ [dd]_{fg}^{}="s_2"\\
\bullet \ar @/^1pc/ [urr]^{p}_{\ }="s" \ar @/_1pc/ [drr]^{\ }="t"_{}&&\\
&&\bullet 
\ar @{=>} @/_.7pc/ "s";"t"_{F'} |>\hole
\ar @{:>} @/^.7pc/ "s";"t"_{K'} |>\hole
\ar @2 "s_2";"t_2"^{\leftact{g}{L}}
}
}
\qquad F' = \hprod{\rightact{f}{G}}{F}, 
 \quad K' = \hprod{\rightact{k}{G}}{K}
$$
and
$$
\raisebox{36pt}{
\xymatrix{
&&\bullet \ar @/^1pc/ [dd]^{k}_{}="t_2"\ar @/_1pc/ [dd]_{f}^{}="s_2"\\
\bullet \ar @/^1pc/ [urr]^{q}_{\ }="s" \ar @/_1pc/ [drr]^{\ }="t"_(.78){}="s_3"   \ar@/_1pc/ [ddrr]_(.8){}^(.8){}="t_3"&&\\
&&\bullet \ar[d]^{g}\\
&&
\ar @{=>} @/_.7pc/ "s";"t"_{F} |>\hole
\ar @{:>} @/^.7pc/ "s";"t"_{K} |>\hole
\ar @2 "s_2";"t_2"^{L}
\ar @{=>} "s_3";"t_3"^{G}
}
}
\equals
\raisebox{36pt}
{
\xymatrix{
&&\bullet \ar @/^1pc/ [dd]^{kg}_{}="t_2"\ar @/_1pc/ [dd]_{fg}^{}="s_2"\\
\bullet \ar @/^1pc/ [urr]^{p}_{\ }="s" \ar @/_1pc/ [drr]^{\ }="t"_{}&&\\
&&\bullet 
\ar @{=>} @/_.7pc/ "s";"t"_{F'} |>\hole
\ar @{:>} @/^.7pc/ "s";"t"_{K'} |>\hole
\ar @2 "s_2";"t_2"^{\rightact{g}{L}}
}
}
\qquad F' = \hprod{\rightact{g}{F}}{G}, 
 \quad K' = \hprod{\rightact{k}{K}}{G}
 \,.
$$

There is one more type of whiskering possible with 2-morphisms,
$$
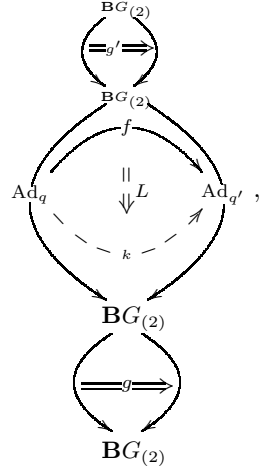
\xymatrix{
  \makebox(8,8){}
  \ar@/^2pc/[rr]|{f}
  \ar@{-->}@/_2pc/[rr]|{\tiny \mbox{$k$}}
  &&
  \makebox(8,8){}
  \ar@{==>}^L (12,3); (12,-3)
  \ar@<+8pt>@/^2.4pc/|<<<<<<<<<<<<{\mathrm{Ad}_{q'}} (12,11.5); (12,-14.5)
  \ar@<-8pt>@/_2.4pc/|<<<<<<<<<<<<{\mathrm{Ad}_{q}} (12,11.5); (12,-14.5)
  \ar@{}|{\small \mbox{$\mathbf{B} G_{(2)}$}} (13,-17.0); (13,-17.0)
  \ar@{}|{\tiny \mbox{$\mathbf{B} G_{(2)}$}} (12,12); (12,12)
  \ar@{}|{\small \mbox{$\mathbf{B} G_{(2)}$}} (13,-17.0)+(0,-18); (13,-17.0)+(0,-18)
  \ar@{}|{\tiny \mbox{$\mathbf{B} G_{(2)}$}} (12,12)+(0,12); (12,12)+(0,12)
  \ar @/_1.2pc/(10,-19.0); (10,-32.0)
  \ar @/^1.2pc/(14,-19.0); (14,-32.0)
  \ar @{=>} (6,-26);(18,-26)|{g}
  \ar @/_.7pc/(9,23); (9,14)
  \ar @/^.7pc/(13,23); (13,14)
  \ar @{=>} (7,19);(15.5,19)|>>>>>>{\tiny \mbox{$g'$}}
}\; ,
$$
which acts in the following way on the components:
$$
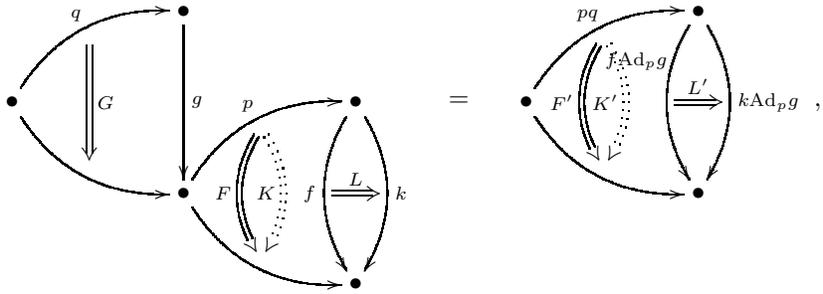
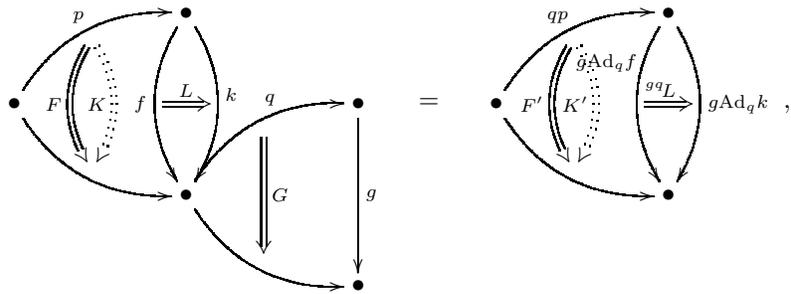
\raisebox{36pt}
{
\xymatrix{
&&\bullet \ar[dd]^{g}                                                                          &&\\
\bullet \ar @/^1pc/ [urr]^{q}_{\ }="s_1" \ar @/_1pc/ [drr]^{\ }="t_1"_{}&& &&
\bullet \ar @/^1pc/ [dd]^{k}_{}="t_3"\ar @/_1pc/ [dd]_{f}^{}="s_3" \\
&&\bullet \ar @/^1pc/ [urr]^{p}_{\ }="s_2" \ar @/_1pc/ [drr]^{\ }="t_2"_{}
\ar @2 "s_1";"t_1"^{G}
\\
&&&&\bullet
\ar @{=>} @/_.7pc/ "s_2";"t_2"_{F}|>\hole
\ar @{:>} @/^.7pc/ "s_2";"t_2"_{K} |>\hole
\ar @2 "s_3";"t_3"^{L}
}
}
\equals
\raisebox{36pt}
{
\xymatrix{
&&\bullet \ar @/^1pc/ [dd]^{k\mathrm{Ad}_pg}_{}="t_2"\ar @/_1pc/ [dd]_(.3){f\mathrm{Ad}_pg\!}^{}="s_2"
\\
\bullet \ar @/^1pc/ [urr]^{pq}_{\ }="s" \ar @/_1pc/ [drr]^{\ }="t"_{}&&
\\
&&\bullet 
\ar @{=>} @/_.7pc/ "s";"t"_{F'} |>\hole
\ar @{:>} @/^.7pc/ "s";"t"_{K'} |>\hole
\ar @2 "s_2";"t_2"^{L'}
}
}\;,
$$
where
\begin{eqnarray*}
F'&=&F\act{fp}{G},\\
K'&=&K\act{kp}{G},\\
L'&=&\act{kp}{G^{-1}}L\act{fp}{G}.
\end{eqnarray*}
and
$$
\raisebox{36pt}
{
\xymatrix{
&&\bullet \ar @/^1pc/ [dd]^(.45){k}_{}="t_3"\ar @/_1pc/ [dd]_{f}^{}="s_3"  &&\\
\bullet \ar @/^1pc/ [urr]^{p}_{\ }="s_1" \ar @/_1pc/ [drr]^{\ }="t_1"_{}&& &&
\bullet \ar[dd]^{g}\\
&&\bullet \ar @/^1pc/ [urr]^(.6){q}_{\ }="s_2" \ar @/_1pc/ [drr]^{\ }="t_2"_{}
\ar @2 "s_2";"t_2"^{G}
\\
&&&&\bullet
\ar @{=>} @/_.7pc/ "s_1";"t_1"_{F}|>\hole
\ar @{:>} @/^.7pc/ "s_1";"t_1"_{K} |>\hole
\ar @2 "s_3";"t_3"^{L}
}
}
\equals
\raisebox{36pt}
{
\xymatrix{
&&\bullet \ar @/^1pc/ [dd]^{g\mathrm{Ad}_qk}_{}="t_2"\ar @/_1pc/ [dd]_(.3){g\mathrm{Ad}_qf\!}^{}="s_2"
\\
\bullet \ar @/^1pc/ [urr]^{qp}_{\ }="s" \ar @/_1pc/ [drr]^{\ }="t"_{}&&
\\
&&\bullet 
\ar @{=>} @/_.7pc/ "s";"t"_{F'} |>\hole
\ar @{:>} @/^.7pc/ "s";"t"_{K'} |>\hole
\ar @2 "s_2";"t_2"^<<<{\act{gq}{L}}
}
}\;,
$$
where
\begin{eqnarray*}
F'&=&G\act{gq}{F},\\
K'&=&G\act{gq}{K}.
\end{eqnarray*}

An important case of this is: 
$$
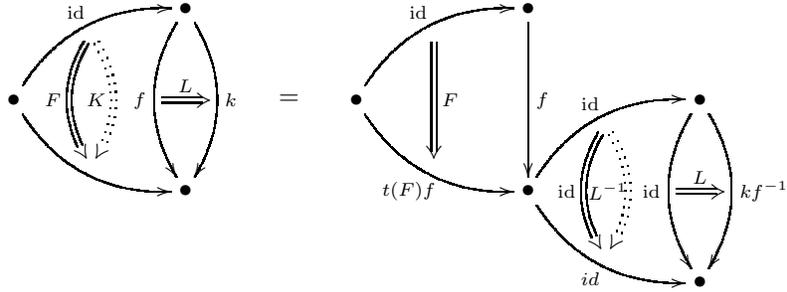
\raisebox{36pt}
{
\xymatrix{
&&\bullet \ar @/^1pc/ [dd]^{k}_{}="t_2"\ar @/_1pc/ [dd]_{f}^{}="s_2"
\\
\bullet \ar @/^1pc/ [urr]^{\id}_{\ }="s" \ar @/_1pc/ [drr]^{\ }="t"_{}&&
\\
&&\bullet 
\ar @{=>} @/_.7pc/ "s";"t"_{F} |>\hole
\ar @{:>} @/^.7pc/ "s";"t"_{K} |>\hole
\ar @2 "s_2";"t_2"^{L}
}
}
\equals
\raisebox{36pt}
{
\xymatrix{
&&\bullet \ar[dd]^{f}                                                                          &&\\
\bullet \ar @/^1pc/ [urr]^{\id}_{\ }="s_1" \ar @/_1pc/ [drr]^{\ }="t_1"_{t(F)f}&& &&
\bullet \ar @/^1pc/ [dd]^{kf^{-1}}_{}="t_3"\ar @/_1pc/ [dd]_{\id}^{}="s_3" \\
&&\bullet \ar @/^1pc/ [urr]^{\id}_{\ }="s_2" \ar @/_1pc/ [drr]^{\ }="t_2"_{id}
\ar @2 "s_1";"t_1"^{F}
\\
&&&&\bullet
\ar @{=>} @/_.7pc/ "s_2";"t_2"_{\id}|>\hole
\ar @{:>} @/^.7pc/ "s_2";"t_2" _{L^{-1}\!\!}|>\hole
\ar @2 "s_3";"t_3"^{L}
}
}.
$$

\subsubsection{Product}

The whiskering along objects
$$
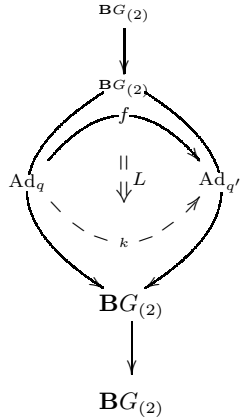
\xymatrix{
  \makebox(8,8){}
  \ar@/^2pc/[rr]|{f}
  \ar@{-->}@/_2pc/[rr]|{\tiny \mbox{$k$}}
  &&
  \makebox(8,8){}
  \ar@{==>}^L (12,3); (12,-3)
  \ar@<+8pt>@/^2.4pc/|<<<<<<<<<<<<{\mathrm{Ad}_{q'}} (12,11.5); (12,-14.5)
  \ar@<-8pt>@/_2.4pc/|<<<<<<<<<<<<{\mathrm{Ad}_{q}} (12,11.5); (12,-14.5)
  \ar@{}|{\small \mbox{$\mathbf{B} G_{(2)}$}} (13,-17.0); (13,-17.0)
  \ar@{}|{\tiny \mbox{$\mathbf{B} G_{(2)}$}} (12,12); (12,12)
  \ar@{}|{\small \mbox{$\mathbf{B} G_{(2)}$}} (13,-17.0)+(0,-13); (13,-17.0)+(0,-13)
  \ar@{}|{\tiny \mbox{$\mathbf{B} G_{(2)}$}} (12,12)+(0,10); (12,12)+(0,10)
  \ar (13,-19.0); (13,-26.0)
  \ar (12,20); (12,14)
}
$$
gives the product of objects with 2-morphisms in the 3-group 
$\mathrm{INN}(G_{(2)})$.
Its action on 2-morphisms, which we have already disucssed, extends 
in a simple way to 3-morphisms:

left whiskering along an object acts as
$$
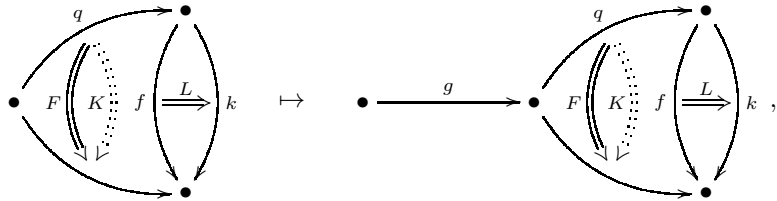
\raisebox{36pt}{
\xymatrix{
&&\bullet \ar @/^1pc/ [dd]^{k}_{}="t_2"\ar @/_1pc/ [dd]_{f}^{}="s_2"
\\
\bullet \ar @/^1pc/ [urr]^{q}_{\ }="s" \ar @/_1pc/ [drr]^{\ }="t"_{}&&
\\
&&\bullet 
\ar @{=>} @/_.7pc/ "s";"t"_{F} |>\hole
\ar @{:>} @/^.7pc/ "s";"t"_{K} |>\hole
\ar @2 "s_2";"t_2"^{L}
}
}
\hspace{10pt}
  \mapsto
\hspace{10pt}
\raisebox{36pt}{
\xymatrix{
&&&&\bullet \ar @/^1pc/ [dd]^{k}_{}="t_2"\ar @/_1pc/ [dd]_{f}^{}="s_2"\\
\bullet \ar[rr]^{g}&&\bullet \ar @/^1pc/ [urr]^{q}_{\ }="s" \ar @/_1pc/ [drr]^{\ }="t"_{}&&\\
&&&&\bullet 
\ar @{=>} @/_.7pc/ "s";"t"_{F} |>\hole
\ar @{:>} @/^.7pc/ "s";"t"_{K} |>\hole
\ar @2 "s_2";"t_2"^{L}
}
}
\,,
$$
while right whiskering along an object
acts as
$$
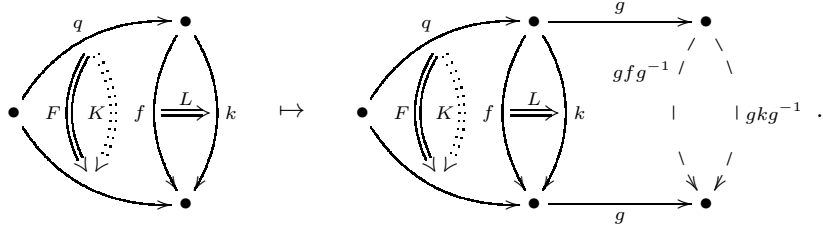
\raisebox{36pt}{
\xymatrix{
&&\bullet \ar @/^1pc/ [dd]^{k}_{}="t_2"\ar @/_1pc/ [dd]_{f}^{}="s_2"
\\
\bullet \ar @/^1pc/ [urr]^{q}_{\ }="s" \ar @/_1pc/ [drr]^{\ }="t"_{}&&
\\
&&\bullet 
\ar @{=>} @/_.7pc/ "s";"t"_{F} |>\hole
\ar @{:>} @/^.7pc/ "s";"t"_{K} |>\hole
\ar @2 "s_2";"t_2"^{L}
}
}
\hspace{10pt}
  \mapsto
\hspace{10pt}
\raisebox{36pt}{
\xymatrix{
&&\bullet \ar @/^1pc/ [dd]^{k}_{}="t_2"\ar @/_1pc/ [dd]_{f}^{}="s_2" \ar[rr]^{g} && \bullet  \ar @{-->}@/^1pc/ [dd]^{gkg^{-1}}\ar@{-->} @/_1pc/ [dd]_(.3){gfg^{-1}}\\
\bullet \ar @/^1pc/ [urr]^{q}_{\ }="s" \ar @/_1pc/ [drr]^{\ }="t"_{}&&&&\\
&&\bullet \ar[rr]_{g} && \bullet
\ar @{=>} @/_.7pc/ "s";"t"_{F} |>\hole
\ar @{:>} @/^.7pc/ "s";"t"_{K} |>\hole
\ar @2 "s_2";"t_2"^{L}
}
}
\,.
$$

To calculate the product of a pair of 2-morphisms, we use the fact that a 2-morphism is uniquely determined by its source and target.

$$
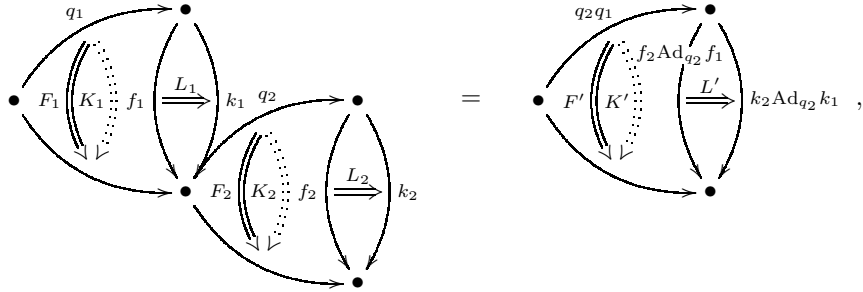
\raisebox{36pt}
{
\xymatrix{
&&\bullet \ar@/^1pc/ [dd]^{k_1}_{}="t_4"\ar @/_1pc/ [dd]_{f_1}^{}="s_4"             &&\\
\bullet \ar @/^1pc/ [urr]^{q_1}_{\ }="s_1" \ar @/_1pc/ [drr]^{\ }="t_1"_{}&& &&
\bullet \ar @/^1pc/ [dd]^{k_2}_{}="t_3"\ar @/_1pc/ [dd]_{f_2}^{}="s_3" \\
&&\bullet \ar @/^1pc/ [urr]^(.6){q_2}_{\ }="s_2" \ar @/_1pc/ [drr]^{\ }="t_2"_{}&&
\ar @{=>} @/_.7pc/ "s_1";"t_1"_{F_1}|>\hole
\ar @{:>} @/^.7pc/ "s_1";"t_1"_{K_1} |>\hole
\ar @2 "s_4";"t_4"^{L_1}
\\
&&&&\bullet
\ar @{=>} @/_.7pc/ "s_2";"t_2"_{F_2}|>\hole
\ar @{:>} @/^.7pc/ "s_2";"t_2"_{K_2} |>\hole
\ar @2 "s_3";"t_3"^{L_2}
}
}
\equals
\raisebox{36pt}
{
\xymatrix{
&&\bullet \ar @/^1pc/ [dd]^{k_2\mathrm{Ad}_{q_2}k_1}_{}="t_2"\ar @/_1pc/ [dd]^{}="s_2"|(.25){f_2\mathrm{Ad}_{q_2}f_1\;\;}
\\
\bullet \ar @/^1pc/ [urr]^{q_2q_1}_{\ }="s" \ar @/_1pc/ [drr]^{\ }="t"_{}&&
\\
&&\bullet 
\ar @{=>} @/_.7pc/ "s";"t"_{F'} |>\hole
\ar @{:>} @/^.7pc/ "s";"t"_{K'} |>\hole
\ar @2 "s_2";"t_2"^{L'}
}
}\;,
$$
with
\begin{eqnarray*}
F' &=& F_2 \act{f_2q_2}{F_1} \\
K' &=& K_2 \act{k_2q_2}{K_1} \\
L' &=&  L_2 \act{f_2q_2}{L_1}
\end{eqnarray*}

\section{Properties of $\mathrm{INN}_0(G_{(2)})$}
\label{properties}

\subsection{Structure morphisms}

We have defined $\mathbf{B}(\mathrm{INN}_0(G_{(2)}))$
essentially as a sub 3-category of $2\mathrm{Cat}$.
The latter is a Gray-category, in that it is a 3-category which is
strict except for the exchange law for composition of 2-morphisms.
Accordingly, also $\mathbf{B}(\mathrm{INN}_0(G_{(2)}))$
is strict except for the exchange law for 2-morphisms.

This means that as a mere 2-groupoid (forgetting the monoidal
structure) $\mathrm{INN}_0(G_{(2)})$ is strict.

\subsubsection{Strictness as a 2-groupoid}

\begin{proposition}
  The underlying 2-groupoid of
  $\mathrm{INN}_0(G_{(2)})$ is strict.
\end{proposition}
\proof
This follows from the rules 
for horizontal and vertical composition of 2-morphisms
in  $\mathrm{INN}_0(G_{(2)})$ -- 
displayed in \ref{composition of 2-morphisms in INNG} --
and the fact that $G_{(2)}$ itself
is a strict 2-group, by assumption.
\endofproof

But the product 2-functor on $\mathrm{INN}_0(G_{(2)})$ 
respects horizontal composition in $\mathrm{INN}_0(G_{(2)})$
only weakly. 
In the language of 2-groups, this corresponds to a failure of the 
Peiffer identity

\subsection{Trivializability}

\begin{proposition}
  \label{connectedness of INNG}
  The 2-groupoid $\mathrm{INN}_0(G_{(2)})$ is connected,
  $$
    \pi_0(\mathrm{INN}_0(G_{(2)})) = 1
    \,.
  $$
\end{proposition}
\proof
  For any two objects $q$ and $q'$ there is the morphism
  $$
 \raisebox{35pt}
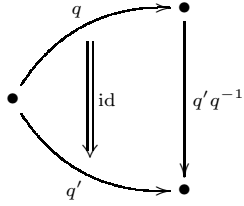
{
\xymatrix{
&&\bullet\ar[dd]^{\gprod{q^{-1}}{q'}}\\
\bullet \ar @/^1pc/ [urr]^{q}_{\ }="s" 
\ar @/_1pc/ [drr]^{\ }="t"_{q'}&&\\
&&\bullet 
\ar @2 "s";"t"^{\id}
}
}       
 $$
\endofproof

\begin{proposition}
  \label{codiscreteness at top level of INNG}
  The Hom-groupoids of the 2-category $\mathrm{INN}_0(G_{(2)})$
  are codiscrete, meaning that they have precisely one morphism
  for every ordered pair of objects.
\end{proposition}
\proof
  By equation (\ref{relation of Inn(G_2) 2-morphisms to source and target})
  there is at most one 2-morphism between any parallel pair of
  morphisms in $\mathrm{INN}_0(G_{(2)})$. For there to be any such 2-morphism
  at all, the two group elements $f$ and $k$ in the diagram above
  (\ref{relation of Inn(G_2) 2-morphisms to source and target}) have
  to satisfy $kf^{-1} \in \mathrm{im}(t)$. But by using the source-target
  matching condition for $F$ and $K$ one readily sees that this is
  always the case.
\endofproof

\begin{theorem}
  The 3-group $\mathrm{INN}_0(G_{(2)})$ is equivalent to the trivial
  3-group. If $G_{(2)}$ is a Lie 2-group, then $\mathrm{INN}_0(G_{(2)})$
  is equivalent to the trivial Lie 3-group even as a Lie 3-group.
\end{theorem}
\proof
  Equivalence of 3-groups $G_{(3)}$, $G'_{(3)}$ is, by definition, 
  that of the corresponding 1-object 3-groupoids
  $\mathbf{B} G_{(3)}$, $\mathbf{B} G'_{(3)}$.
  For showing equivalence with the trivial 3-group, it suffices to 
  exhibit a pseudonatural transformation of 3-functors
  $$
    \id_{\mathbf{B} (\mathrm{INN}_0(G_{(2)}))} 
    \to
    I_{\mathbf{B} (\mathrm{INN}_0(G_{(2)}))}
    \,,
  $$
  where $I_{\mathbf{B} (\mathrm{INN}_0(G_{(2)}))}$ sends everything to the
  identity on the single object of $\mathbf{B} \mathrm{INN}_0(G_{(2)})$.
  Such a transformation is obtained by sending the single object to the
  identity 1-morphism on that object and sending any 1-morphism $q$ to the
  2-morphism $q \to \id$ from prop \ref{connectedness of INNG}.
  By prop \ref{codiscreteness at top level of INNG} this implies 
  the existence of a
  unique assignment of a 3-morphism to any 2-morphism such that we do 
  indeed obtain the component map of a pseudonatural transformation
  of 3-functors. By construction, this is clearly smooth when $G_{(2)}$
  is Lie. 
\endofproof

\subsection{Universality}

\begin{theorem}
  We have a short exact sequence of strict 2-groupoids
  $$
    \xymatrix{
      G_{(2)}
      \ar@{^{(}->}[r]
      &
      \mathrm{INN}_0(G_{(2)})
      \ar@{->>}[r]
      &
      \mathbf{B} G_{(2)}
    }
      \,.
  $$
\end{theorem}
\proof
This is just proposition \ref{exact sequence}, after noticing that
$$
  \mathrm{Mor}(\mathbf{B} G_{(2)})
  =
  G_{(2)}
  \,.
$$
\endofproof

So the strict inclusion 2-functor on the left is
$$
  \left(
  \xymatrix{
    g
    \ar[r]^h
    &
    g'
  }
  \right)
  \hspace{7pt}
    \mapsto
  \hspace{7pt}
\raisebox{36pt}{
\xymatrix{
&&\bullet \ar @/^1pc/ [dd]^{\id}_{}="t_2"\ar @/_1pc/ [dd]|{\id}^{}="s_2"
\\
\bullet \ar @/^1pc/ [urr]^{g}_{\ }="s" \ar @/_1pc/ [drr]^{\ }="t"_{g'}&&
\\
&&\bullet 
\ar @{=>} @/_.7pc/ "s";"t"_{h} |>\hole
\ar @{:>} @/^.7pc/ "s";"t"_{h} |>\hole
\ar @2 "s_2";"t_2"^{\id}
}
}
\,,
$$
while the strict surjection 2-functor on the right is
$$
\raisebox{36pt}{
\xymatrix{
&&\bullet \ar @/^1pc/ [dd]^{k}_{}="t_2"\ar @/_1pc/ [dd]_{f}^{}="s_2"
\\
\bullet \ar @/^1pc/ [urr]^{q}_{\ }="s" \ar @/_1pc/ [drr]^{\ }="t"_{}&&
\\
&&\bullet 
\ar @{=>} @/_.7pc/ "s";"t"_{F} |>\hole
\ar @{:>} @/^.7pc/ "s";"t"_{K} |>\hole
\ar @2 "s_2";"t_2"^{L}
}
}
\hspace{7pt}
  \mapsto
\hspace{7pt}
\xymatrix{
  \bullet
  \ar@/^1.8pc/[rr]^{f}_{\ }="s"
  \ar@/_1.8pc/[rr]_{k}^{\ }="t"
  &&
  \bullet
  \ar@{=>}^L "s"; "t"  
}
\,.
$$

\subsection{The corresponding 2-crossed module}
\label{The corresponding 2-crossed module}

We now extract the structure of a 2-crossed module from 
$\mathrm{INN}_0(G_{(2)})$. First, let
$$
\mathrm{Mor}_1^{I}  = \mathrm{Mor}_1(\mathrm{INN}_0(G_{(2)}))|_{s^{-1}
(\mathrm{Id})}
$$
and
$$
\mathrm{Mor}_2^{I} =   \mathrm{Mor}_2(\mathrm{INN}_0(G_{(2)}))|_{s^{-1}
(\id_{\mathrm{Id}})}
$$
be subgroups of the 1- and 2-morphisms of $\mathrm{INN}_0(G_{(2)})$
respectively.

\begin{proposition}\label{1-morphisms from id}
The group of 1-morphisms in $\mathrm{INN}_0(G_{(2)})$ starting at the
identity object form the semidirect product group
$$
  \mathrm{Mor}_1^{I}  = G \ltimes H
$$
under the identitfication
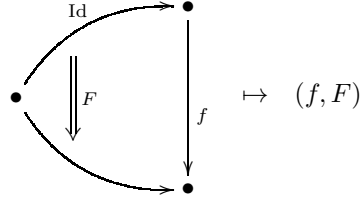
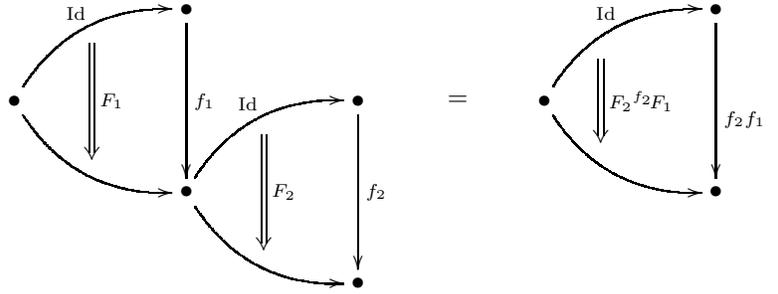
$$
\raisebox{36pt}
{
\xymatrix{
&&
\bullet
\ar[dd]^(.6){f}
\\
\bullet
 \ar @/^1pc/ [urr]^{\mathrm{Id}}_(.35){\ }="s"
 \ar @/_1pc/ [drr]_{}^(.35){\ }="t"_{}&&\\
&&\bullet
\ar@{=>} "s";"t"^{F}
}
}
\hspace{7pt}
  \mapsto
\hspace{7pt}
  (f,F)
$$
in that
\end{proposition}
\begin{eqnarray*}
\raisebox{36pt}
{
\xymatrix{
&&\bullet \ar[dd]^
{f_1}
       &&\\
\bullet
\ar @/^1pc/ [urr]^{\mathrm{Id}}_{\ }="s_1"
\ar @/_1pc/ [drr]_{}^{\ }="t_1"_{}
&&
&&\bullet\ar[dd]^{f_2}\\
&&\bullet
\ar @/^1pc/ [urr]^{\mathrm{Id}}_{\ }="s_2"
\ar @/_1pc/ [drr]_{}^{\ }="t_2"_{}
\ar @2 "s_1";"t_1"^{F_1}
\\
&&&&\bullet
\ar @2 "s_2";"t_2"^{F_2}
}
}
&\equals&
\raisebox{36pt}
{
\xymatrix{
&&
\bullet
\ar[dd]^(.6){f_2f_1}
\\
\bullet
 \ar @/^1pc/ [urr]^{\mathrm{Id}}_(.35){\ }="s"
 \ar @/_1pc/ [drr]_{}^(.35){\ }="t"_{}&&\\
&&\bullet
\ar@{=>} "s";"t"^{F_2\act{f_2}{F_1}}
}
}
\,.
\end{eqnarray*}
\proof
  Use composition in $\mathbf{B} G_{(2)}$.
\endofproof

We have the obvious group homomorphism which is just the restriction
of the target map
$$
  \partial_1 : \mathrm{Mor}_1^{I}  \to \mathrm{Obj} := \mathrm{Obj}
(\mathrm{INN}_0(G_{(2)}))
$$
given by
$$
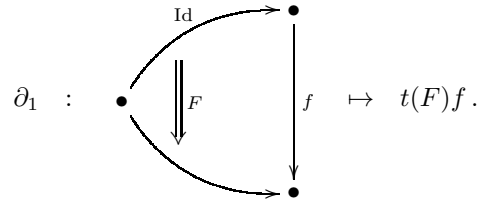
\partial_1
\hspace{7pt}
  :
\hspace{7pt}
\raisebox{36pt}
{
\xymatrix{
&&
\bullet
\ar[dd]^{f}
\\
\bullet
 \ar @/^1pc/ [urr]^{\mathrm{Id}}_(.35){\ }="s"
 \ar @/_1pc/ [drr]_{}^(.35){\ }="t"_{}&&\\
&&\bullet
\ar@{=>} "s";"t"^{F}
}
}
\hspace{7pt}
  \mapsto
\hspace{7pt}
  t(F)f
  \,.
$$
This and the following constructions are to be compared with
definition \ref{nonabelian mapping cone}. There is an obvious action
on $\mathrm{Mor}_1^{I}$:

\begin{equation}
\label{action of objects on arrows}
\hspace{-2cm}
\raisebox{36pt}
{
\xymatrix{
&&
 \bullet\ar[dd]^{f}
  \\
  \bullet
  \ar @/^1pc/ [urr]^{\id}_{\ }="s" \ar @/_1pc/ [drr]^{\ }="t"_{}&&\\
&&\bullet
\ar @2 "s";"t"^{F}
}
}
\hspace{8pt}
  \mapsto
\hspace{8pt}
\raisebox{36pt}
{
\xymatrix{
&&&&\bullet \ar[dd]^{f} \ar[rr]^{g} \ar@{}[ddrr] |{=} && \bullet
 \ar@{-->}[dd]^{\gprod{g^{-1}}{\gprod{f}{g}}}\\
\bullet \ar[rr]^{g^{-1}} &&\bullet
 \ar @/^1pc/ [urr]^{\id}_{\ }="s"
\ar @/_1pc/ [drr]^{\ }="t"_{}&&&&\\
&&&&
 \bullet \ar[rr]_g &&\bullet
\ar @2 "s";"t"^{F}
}
}
\equals
\raisebox{36pt}
{
\xymatrix{
&&\bullet\ar[dd]^{\gprod{g^{-1}}{\gprod{f}{g}}}\\
\bullet
  \ar @/^1pc/ [urr]^{\id}_{\ }="s"
 \ar @/_1pc/ [drr]^{\ }="t"_{}&&\\
&&\bullet
\ar @2 "s";"t"^{\rightact{g}{F}}
}
}
\end{equation}
This action almost gives us a crossed module $\mathrm{Mor}_1^{I} \to
\mathrm{Obj}$.
But not quite, since the Peiffer identity holds only up to 3-
isomorphism.

To see this, let

$$
g
\equals
\del_1\left(
\raisebox{36pt}
{
\xymatrix{
&&\bullet\ar[dd]^{h}\\
\bullet \ar @/^1pc/ [urr]^{\id}_{\ }="s" \ar @/_1pc/ [drr]^{\ }="t"_
{g}&&\\
&&\bullet
\ar @2 "s";"t"^{H}
}
}
\right) \equals t(H)h.
$$
For the Peiffer identity to hold we need the action (\ref{action of
objects on arrows}) to be equal to the adjoint action of the 2-cell $
(h,H;id)$. To see that this fails, first notice that
the inverse of the approriate 2-cell considered as an element in the
group $\mathrm{Mor}_1^{I}$ is
$$
\left(
\raisebox{36pt}
{
\xymatrix{
&&\bullet\ar[dd]^{h}\\
\bullet \ar @/^1pc/ [urr]^{\id}_{\ }="s" \ar @/_1pc/ [drr]^{\ }="t"_(.6){t(H)h}&&\\
&&\bullet
\ar @2 "s";"t"^{H}
}
}
\right)^{-1}
\equals
\raisebox{36pt}
{
\xymatrix{
&&\bullet\ar[dd]^(.35){h^{-1}}\\
\bullet \ar @/^1pc/ [urr]^{\id}_{\ }="s" \ar @/_1pc/ [drr]^{\ }="t"_(.7){h^{-1}t(H)^{-1}}&&\\
&&\bullet
\ar @2 "s";"t"^{{}^{h^{-1}}H^{-1}}
}
}
\,.
$$
Therefore the conjugation is

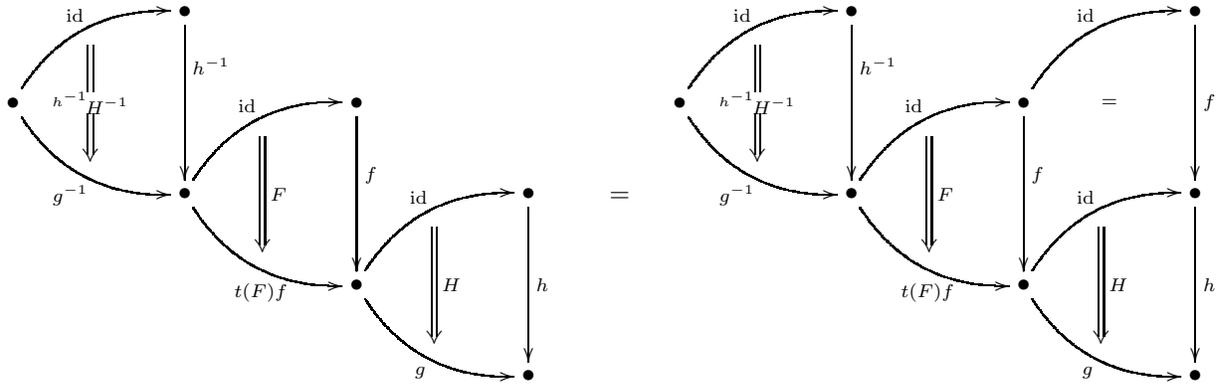
\begin{eqnarray}
\hspace{-2cm}
\raisebox{72pt}
{
\xymatrix{
&&\bullet \ar[dd]^(.3){h^
{-1}}
&&     &&\\
\bullet \ar @/^1pc/ [urr]^{\id}_{\ }="s_1" \ar @/_1pc/ [drr]^{\ }
="t_1"_{g^{-1}}&& &&\bullet\ar[dd]^(.4){f}  &&\\
&&\bullet \ar @/^1pc/ [urr]^{\id}_{\ }="s_2" \ar @/_1pc/ [drr]^{\ }
="t_2"_(.6){t(F)f}&& && \bullet \ar[dd]^{h}
\ar @2 "s_1";"t_1"|{\act{h^{-1}}{H^{-1}}}
\\
&&&&\bullet  \ar @/^1pc/ [urr]^{\id}_{\ }="s_3" \ar @/_1pc/ [drr]^{\ }
="t_3"_{g} &&\\
&&&&&& \bullet
\ar @2 "s_2";"t_2"^{F}
\ar @2 "s_3";"t_3"^{H}
}
}
&\hspace{10pt}=&
\raisebox{72pt}
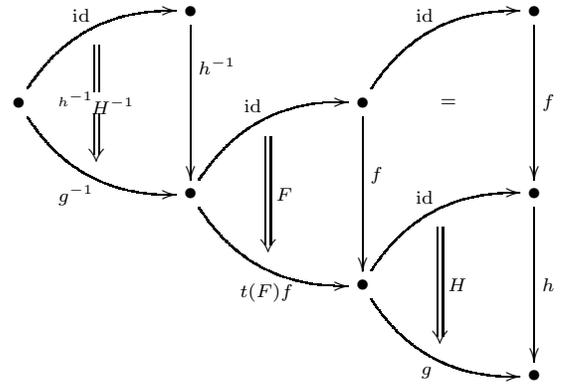
{
\xymatrix{
&&\bullet \ar[dd]^(.3){h^
{-1}}
&&     &&\bullet \ar[dd]^{f}\\
\bullet \ar @/^1pc/ [urr]^{\id}_{\ }="s_1" \ar @/_1pc/ [drr]^{\ }
="t_1"_{g^{-1}}&& &&\bullet\ar[dd]^(.4){f}  \ar @/^1pc/ [urr]^{\id}
\ar@{}[rr] |{=} &&\\
&&\bullet \ar @/^1pc/ [urr]^{\id}_{\ }="s_2" \ar @/_1pc/ [drr]^{\ }
="t_2"_(.6){t(F)f}&& && \bullet \ar[dd]^{h}
\ar @2 "s_1";"t_1"|{\act{h^{-1}}{H^{-1}}}
\\
&&&&\bullet  \ar @/^1pc/ [urr]^{\id}_{\ }="s_3" \ar @/_1pc/ [drr]^{\ }
="t_3"_{g} &&\\
&&&&&& \bullet
\ar @2 "s_2";"t_2"^{F}
\ar @2 "s_3";"t_3"^{H}
}
}
\nonumber\\ %End first line
&&\nonumber\\ %Gap
&&\nonumber\\ %Gap
&\hspace{10pt}=&
\raisebox{72pt}
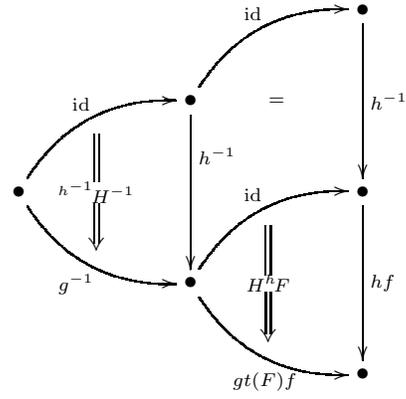
{
\xymatrix{
&&&&\bullet\ar[dd]^{h^{-1}}\\
&&\bullet \ar[dd]^(.3){h^{-1}}   \ar @/^1pc/ [urr]^{\id}  \ar@{}[rr] |
{=}                                                     &&\\
\bullet  \ar @/^1pc/ [urr]^{\id}_{\ }="s_1" \ar @/_1pc/ [drr]^{\ }
="t_1"_{g^{-1}}&& &&\bullet\ar[dd]^{\gprod{f}{h}}\\
&&\bullet \ar @/^1pc/ [urr]^{\id}_{\ }="s_2" \ar @/_1pc/ [drr]^{\ }
="t_2"_(.6){\gprod{t(F)f}{g}}
\ar @2 "s_1";"t_1"|{\act{h^{-1}}{H^{-1}}}
\\
&&&&\bullet
\ar @2 "s_2";"t_2"|{
                             \hprod{\rightact{h}{F}}{H}
                             }
}
}
\nonumber\\ %End second line
&&\nonumber\\ %Gap
&&\nonumber\\ %Gap
&\hspace{10pt}=&
\raisebox{36pt}
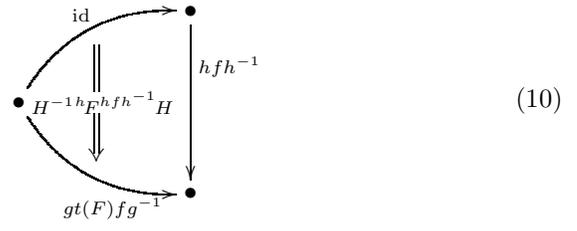
{
\xymatrix{
&&\bullet\ar[dd]^(.3){hfh^{-1}}\\
\bullet \ar @/^1pc/ [urr]^{\id}_{\ }="s" \ar @/_1pc/ [drr]^{\ }="t"_(.7){gt(F)fg^{-1}}&&\\
&&\bullet
\ar @2 "s";"t"|{\;\;
                        \hprod{\rightact{hfh^{-1}}{H}}{\hprod
{\rightact{h}{F}}{H^{-1}}}
                        }
}
}
\label{conjugation action}
\\ %End third line
&&\nonumber\\ %Gap
&&\nonumber\\ %Gap
&\hspace{10pt}\neq&
\raisebox{36pt}
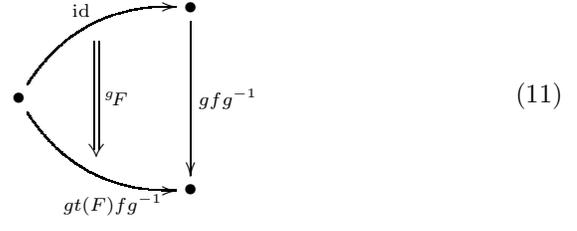
{
\xymatrix{
&&\bullet\ar[dd]^{gfg^{-1}}\\
\bullet \ar @/^1pc/ [urr]^{\id}_{\ }="s" \ar @/_1pc/ [drr]^{\ }="t"_(.7){gt(F)fg^{-1}}&&\\
&&\bullet
\ar @2 "s";"t"^{\act{g}{F}}
}
}
\label{action}
\end{eqnarray}
Though the Peiffer identity does not hold, both actions give rise to
2-cells with the same source and target, and hence define a 3-cell $P
$. Denote the 2-cell (\ref{conjugation action}) by $(c,C;id)$ and the
2-cell (\ref{action}) by $(a,A;id)$ (for \emph{c}onjugation and \emph
{a}ction respectively).
$$
\xymatrix{
&&\bullet \ar @/^1pc/ [dd]^{a}_{}="t_2"\ar @/_1pc/ [dd]_{c}^{}="s_2"\\
\bullet \ar @/^1pc/ [urr]^{\id}_{\ }="s" \ar @/_1pc/ [drr]^{\ }="t"_{}
&&\\
&&\bullet
\ar @{=>} @/_.7pc/ "s";"t"_{C} |>\hole
\ar @{:>} @/^.7pc/ "s";"t"_{A} |>\hole
\ar @2 "s_2";"t_2"^{P}
}
$$
Then
\begin{eqnarray*}
P &=& A^{-1} C\\
&=&\act{g}{F^{-1}}\left(H\rightact{h}{F}\rightact{hfh^{-1}}{H^{-1}}
\right) \\
&=&\act{g}{F^{-1}}\left(\rightact{t(H)h}{F}H\rightact{hfh^{-1}}{H^
{-1}}\right)\\
&=&\act{g}{F^{-1}}\left(\rightact{g}{F}H\rightact{hfh^{-1}}{H^{-1}}
\right)\\
&=&H\rightact{hfh^{-1}}{H^{-1}}
\end{eqnarray*}
However, what we really want is the Peiffer lifting, which will be a
3-cell with source the identity 2-cell. Hence,

\begin{proposition}\label{2-morphisms from id}
The group of 2-morphisms in $\mathrm{INN}_0(G_{(2)})$ starting at the
identity arrow on the identity object form the group
$$
\mathrm{Mor}_2^{I}  = H
$$
under the identitfication
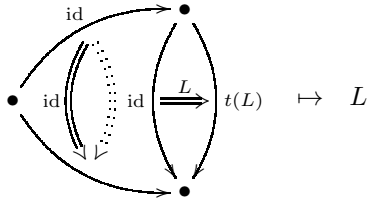
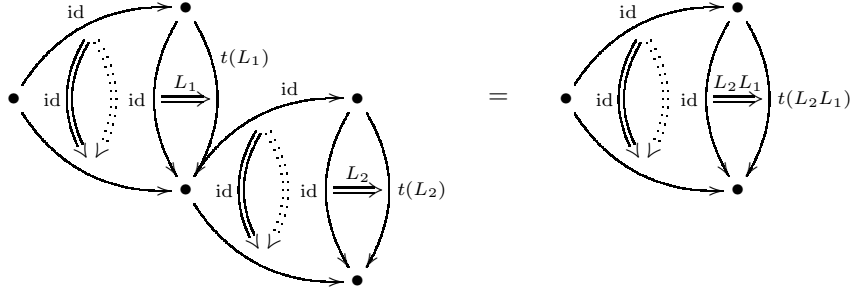
$$
\raisebox{36pt}
{
\xymatrix{
&&\bullet \ar @/^1pc/ [dd]^{t(L)}_{}="t_2"\ar @/_1pc/ [dd]_{\id}^{}
="s_2"
\\
\bullet \ar @/^1pc/ [urr]^{\id}_{\ }="s" \ar @/_1pc/ [drr]^{\ }="t"_{}&&
\\
&&\bullet
\ar @{=>} @/_.7pc/ "s";"t"_{\id} |>\hole
\ar @{:>} @/^.7pc/ "s";"t"_{} |>\hole
\ar @2 "s_2";"t_2"^{L}
}}
\hspace{7pt}
  \mapsto
\hspace{7pt}
  L
$$
in that
\end{proposition}
$$
\raisebox{36pt}
{
\xymatrix{
&&\bullet \ar @/^1pc/ [dd]^(.3){t(L_1)}_{}="t_2"\ar @/_1pc/ [dd]_{\id}
^{}="s_2"&&
\\
\bullet \ar @/^1pc/ [urr]^{\id}_{\ }="s" \ar @/_1pc/ [drr]^{\ }="t"_{}
&&&&\bullet \ar @/^1pc/ [dd]^{t(L_2)}_{}="t_3"\ar @/_1pc/ [dd]_{\id}^
{}="s_3"
\\
&&\bullet  \ar @/^1pc/ [urr]^(.7){\id}_{\ }="s_4" \ar @/_1pc/ [drr]^
{\ }="t_4"_{}&&
\\
&&&&\bullet
\ar @{=>} @/_.7pc/ "s";"t"_{\id} |>\hole
\ar @{:>} @/^.7pc/ "s";"t"_{} |>\hole
\ar @2 "s_2";"t_2"^{L_1}
\ar @{=>} @/_.7pc/ "s_4";"t_4"_{\id} |>\hole
\ar @{:>} @/^.7pc/ "s_4";"t_4"_{} |>\hole
\ar @2 "s_3";"t_3"^{L_2}
}
}
\equals
\raisebox{36pt}
{
\xymatrix{
&&\bullet \ar @/^1pc/ [dd]^{t(L_2L_1)}_{}="t_2"\ar @/_1pc/ [dd]_{\id}^
{}="s_2"
\\
\bullet \ar @/^1pc/ [urr]^{\id}_{\ }="s" \ar @/_1pc/ [drr]^{\ }="t"_{}&&
\\
&&\bullet
\ar @{=>} @/_.7pc/ "s";"t"_{\id} |>\hole
\ar @{:>} @/^.7pc/ "s";"t"_{} |>\hole
\ar @2 "s_2";"t_2"^{L_2L_1}
}
}
$$
\proof
  Use the multiplication of 2-morphisms.
\endofproof

 So, we whisker the 3-cell $(P;a,A;\id)$ above with the inverse of $
(a,A;\id)$:

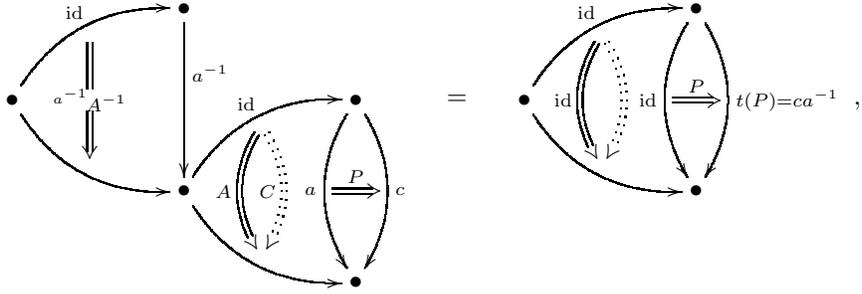
$$
\raisebox{36pt}
{
\xymatrix{
&&\bullet \ar[dd]^(.35){a^
{-1}}
       &&\\
\bullet \ar @/^1pc/ [urr]^{\id}_{\ }="s_1" \ar @/_1pc/ [drr]^{\ }
="t_1"_{}&& &&
\bullet \ar @/^1pc/ [dd]^{c}_{}="t_3"\ar @/_1pc/ [dd]_{a}^{}="s_3" \\
&&\bullet \ar @/^1pc/ [urr]^{\id}_{\ }="s_2" \ar @/_1pc/ [drr]^{\ }
="t_2"_{}
\ar @2 "s_1";"t_1"|{
                               \act{a^{-1}}{A^{-1}}
                               }
\\
&&&&\bullet
\ar @{=>} @/_.7pc/ "s_2";"t_2"_{A}|>\hole
\ar @{:>} @/^.7pc/ "s_2";"t_2"_{C} |>\hole
\ar @2 "s_3";"t_3"^{P}
}
}
\equals
\raisebox{36pt}
{
\xymatrix{
&&\bullet \ar @/^1pc/ [dd]^{t(P)=ca^{-1}}_{}="t_2"\ar @/_1pc/ [dd]_
{\id}^{}="s_2"
\\
\bullet \ar @/^1pc/ [urr]^{\id}_{\ }="s" \ar @/_1pc/ [drr]^{\ }="t"_{}&&
\\
&&\bullet
\ar @{=>} @/_.7pc/ "s";"t"_{\id} |>\hole
\ar @{:>} @/^.7pc/ "s";"t"_{}|>\hole
\ar @2 "s_2";"t_2"^{P}
}
}\;,
$$
and the back face is necessarily $P^{-1}$.

\begin{definition}[Peiffer lifting]
Define the map
$$
  \{\cdot, \cdot\} : \mathrm{Mor}_1^{I} \times \mathrm{Mor}_1^{I}
\to \mathrm{Mor}_2^{I}
$$
by
$$
\left\lbrace
\raisebox{36pt}
{
\xymatrix{
&&\bullet\ar[dd]^{h}\\
\bullet \ar @/^1pc/ [urr]^{\id}_{\ }="s" \ar @/_1pc/ [drr]^{\ }="t"_(.6){}&&\\
&&\bullet
\ar @2 "s";"t"^{H}
}
},
\raisebox{36pt}
{
\xymatrix{
&&\bullet\ar[dd]^{f}\\
\bullet \ar @/^1pc/ [urr]^{\id}_{\ }="s" \ar @/_1pc/ [drr]^{\ }="t"_(.7){}&&\\
&&\bullet
\ar @2 "s";"t"^{F}
}
}
\right\rbrace
\equals
\raisebox{36pt}
{
\xymatrix{
&&\bullet \ar @/^1pc/ [dd]^{t(P)}_{}="t_2"\ar @/_1pc/ [dd]_{\id}^{}
="s_2"
\\
\bullet \ar @/^1pc/ [urr]^{\id}_{\ }="s" \ar @/_1pc/ [drr]^{\ }="t"_{}&&
\\
&&\bullet
\ar @{=>} @/_.7pc/ "s";"t"_{\id} |>\hole
\ar @{:>} @/^.7pc/ "s";"t"_{}|>\hole
\ar @2 "s_2";"t_2"^{P}
}
}\;,\qquad P = H\rightact{hfh^{-1}}{H^{-1}}
$$
\end{definition}

Now define the homomorphism
$$
  \del_2 : \mathrm{Mor}_2^{I} \to \mathrm{Mor}_1^{I}
$$
by
$$
  \del_2
  :
  \raisebox{36pt}
{
\xymatrix{
&&\bullet \ar @/^1pc/ [dd]^{t(L)}_{}="t_2"\ar @/_1pc/ [dd]_{\id}^{}
="s_2"
\\
\bullet \ar @/^1pc/ [urr]^{\id}_{\ }="s" \ar @/_1pc/ [drr]^{\ }="t"_{}&&
\\
&&\bullet
\ar @{=>} @/_.7pc/ "s";"t"_{\id} |>\hole
\ar @{:>} @/^.7pc/ "s";"t"_{} |>\hole
\ar @2 "s_2";"t_2"^{L}
}}
  \hspace{7pt}
  \mapsto
  \hspace{7pt}
\raisebox{36pt}{
\xymatrix{
&&\bullet\ar[dd]^{t(L)}\\
\bullet
\ar @/^1pc/ [urr]^{\mathrm{Id}}_{\ }="s"
\ar @/_1pc/ [drr]^{\ }="t"_{\mathrm{Id}}&&\\
&&\bullet
\ar @2 "s";"t"^{L^{-1}}
}
}\;,
$$
which is again the restriction of the target map. Note there is an
action of $\mathrm{Obj}$ on $\mathrm{Mor}_2^{I}$:
\begin{equation}
\label{action of arrows on 2-arrows}
\hspace{-3cm}
\raisebox{36pt}{
\xymatrix{
&&\bullet \ar @/^1pc/ [dd]^{t(L)}_{}="t_2"\ar @/_1pc/ [dd]_{\id}^{}
="s_2"
\\
\bullet \ar @/^1pc/ [urr]^{\id}_{\ }="s" \ar @/_1pc/ [drr]^{\ }="t"_{}&&
\\
&&\bullet
\ar @{=>} @/_.7pc/ "s";"t"_{\id} |>\hole
\ar @{:>} @/^.7pc/ "s";"t"_{} |>\hole
\ar @2 "s_2";"t_2"^{L}
}
}
\hspace{10pt}
  \mapsto
\hspace{10pt}
\raisebox{36pt}{
\xymatrix{
&&&&\bullet \ar @/^1pc/ [dd]^{t(L)}_{}="t_2"\ar @/_1pc/ [dd]_{\id}^{}
="s_2" \ar[rr]^{g} &&
\bullet  \ar @{-->}@/^1pc/ [dd]^{gt(L)g^{-1}}\ar@{-->} @/_1pc/ [dd]_(.3){\id}
\\
\bullet \ar[rr]^{g}&&\bullet \ar @/^1pc/ [urr]^{\id}_{\ }="s" \ar @/
_1pc/ [drr]^{\ }="t"_{}&&\\
&&&&\bullet \ar[rr]_{g} && \bullet
\ar @{=>} @/_.7pc/ "s";"t"_{\id} |>\hole
\ar @{:>} @/^.7pc/ "s";"t"_{} |>\hole
\ar @2 "s_2";"t_2"^{L}
}
}
\equals
\raisebox{36pt}{
\xymatrix{
&&\bullet \ar @/^1pc/ [dd]^{t(\act{g}{L})}_{}="t_2"\ar @/_1pc/ [dd]_
{\id}^{}="s_2"
\\
\bullet \ar @/^1pc/ [urr]^{\id}_{\ }="s" \ar @/_1pc/ [drr]^{\ }="t"_{}&&
\\
&&\bullet
\ar @{=>} @/_.7pc/ "s";"t"_{\id} |>\hole
\ar @{:>} @/^.7pc/ "s";"t"_{} |>\hole
\ar @2 "s_2";"t_2"^{\act{g}{L}}
}
}
\end{equation}

Clearly $\del_2 \circ \del_1$ is the constant map at the identity,
and $\mathrm{im}\ \del_2$ is a normal subgroup of $\mathrm{ker}\
\del_1$, so
\begin{eqnarray}
\label{2-crossed module of INNG2}
\xymatrix{
\mathrm{Mor}_2^{I} \ar[r]^{\del_2} & \mathrm{Mor}_1^{I} \ar[r]^
{\del_1} & \mathrm{Obj}
}
\end{eqnarray}
is a sequence. We let the action of $\mathrm{Obj}$ on the other two
groups be as described above in (\ref{action of objects on arrows})
and (\ref{action of arrows on 2-arrows}), and the maps $\del_2$ and $
\del_1$ are clearly equivariant for this action.

\begin{proposition}
  The map $\{\cdot,\cdot\}$ does indeed satisfy the properties of a
  Peiffer lifting, and (\ref{2-crossed module of INNG2}) is a 2-
crossed module.
\end{proposition}
\proof
  The first condition holds by definition, the second and the last
  one are easy to check. The others are tedious. It is easy, using
the crossed module properties of $H \to G$, to calculate that the
actions of $\mathrm{Mor}_1^{I}$ on $\mathrm{Mor}_2^{I}$ as defined
from $\mathrm{INN}_0(G_{(2)})$ and as defined via $\{\cdot,\cdot\}$ are
the same.
\endofproof

Since $\mathrm{im}\ \del_2 = \mathrm{ker}\ \del_1$, $\del_2$ is
injective and $\del_1$ is onto, this shows that (\ref{2-crossed
module of INNG2}) has trivial homology and provides us with another
proof that $\mathrm{INN}_0(G_{(2)})$ is contractible.

\subsubsection{Relation to the mapping cone of $H \to G$}

Given a crossed square
  $$
    \xymatrix{
      L
       \ar[rr]^{f}
       \ar[dd]_{u}
       &&
       M
       \ar[dd]^{v}
      \\
      \\
      N \ar[rr]^{g}&& P
    }
  $$
with structure map $h: N \times M \to L$, Conduch\'e 
\cite{Conduche}gives the
Peiffer lifting of the mapping cone
$$
\xymatrix{
L \ar[r] &N \ltimes M \ar[r] & P
}
$$
 as
$$
\{(g,h),(k,l)\} = h(gkg^{-1},h).
$$
Recall from \ref{Mapping cones of crossed modules}
that the identity map on $t:H \to G$ is a crossed square with
$$
h(g,h) = h\act{g}{h^{-1}},
$$
so the mapping cone is a 2-crossed module
$$
\xymatrix{
H \ar[r]^<<<<<<{\partial_2} & G \ltimes H \ar[r]^<<<<<<{\partial_1} & G,
}
$$
where
$$d_2(h) = (t(h),h^{-1}), \qquad d_1(g,h) = t(h)g,$$
and with Peiffer lifting
$$
\{(g_1,h_1),(g_2,h_2)\} = h_1\act{g_1g_2g_1^{-1}}{h_1^{-1}}.
$$
which is what we found for $\mathrm{INN}_0(G_{(2)})$.

More precisely,

\begin{definition}
A morphism $\psi$ of 2-crossed modules is a map of the underlying
complexes
$$
\xymatrix{
L_1 \ar[r]^{\del_2}  \ar[d]_{\psi_L} & M_1 \ar[r]^{\del_1}  \ar[d]_
{\psi_M} & N_1 \ar[d]_{\psi_N} \\
L_2 \ar[r]_{\partial_2} & M_2 \ar[r]_{\partial_1} & N_2
}
$$
such that $\psi_L,\ \psi_M$ and $\psi_N$ are equivariant for the $N$-
and $M$-actions, and
$$
\{\psi_M(\cdot),\psi_M(\cdot)\}_2 = \psi_L(\{\cdot,\cdot\}_1).
$$
\end{definition}

Using propositions \ref{1-morphisms from id} and \ref{2-morphisms
from id}, we have a map

$$
\xymatrix{
\mathrm{Mor}_2^{I} \ar[r]^{\del_2} \ar[d]_{\simeq} & \mathrm{Mor}_1^
{I} \ar[r]^{\del_1} \ar[d]^{\simeq} & \mathrm{Obj} \ar[d]^{\simeq} \\
H \ar[r]_{d_2} & G \ltimes H \ar[r]_{d_1} & G
}
$$

and the actions and Peiffer lifting agree, so

\begin{proposition}
The 2-crossed module associated to $\mathrm{INN}_0(G_{(2)})$ is
isomorphic to the mapping cone of the identity map on the crossed
module associated to $G_{(2)}$.
\end{proposition}

\section{Universal $n$-bundles}
\label{universal n-bundles}

In order to put the relevance of the 3-group $\mathrm{INN}_0(G_{(2)})$
in perspective, we further illuminate our statement, 
\ref{the exact sequence}, 
that $\mathrm{INN}_0(G_{(2)})$ 
\emph{plays the role of the universal $G_{(2)}$-bundle}. 
An exhaustive discussion will be given elsewhere.

\subsection{Universal 1-bundles in terms of $\mathrm{INN}(G)$}

Let $\pi : Y \to X$ be a good cover of a space $X$ and write
$Y^{[2]} := Y \times_X Y$ for the corresponding groupoid.

\begin{definition}[$G$-cocycles]
  A $G$-(1-)cocycle on $X$ is a functor
  $$
    g : Y^{[2]} \to \mathbf{B} G
    \,.
  $$
\end{definition}
This functor can be understood as arising from a choice 
$$
  \xymatrix{
    \pi^* P \ar[r]^t_\sim & Y \times G
  }
$$
of trivialization of a principal right $G$-bundle $P \to X$
(which is essentially just a map to $G$) as
$$
  g := \pi_2^* t \circ \pi_1^* t^{-1}
  \,,
$$
by noticing that $G$-equivariant isomorphisms
$$
  G \to G
$$
are in bijection with elements of $G$
$$
  g(x,y) : h \mapsto g(x,y) h
$$
acting from the left.

\begin{observation}[$G$-bundles as morphisms of sequences of groupoids]
  Given a $G$-cocycle on $X$ as above, its pullback along the exact sequence
  $$
    \xymatrix{
      G \ar[r] & \mathrm{INN}(G) \ar[r] & \mathbf{B} G
    }
    \,,
  $$
  which we write as
$$
  \raisebox{40pt}{
  \xymatrix{
    Y \times G
    \ar[rr]
    \ar[dd]
    &&
    Y^{[2]} \times_{g} \mathrm{INN}(G)
    \ar[rr]
    \ar[dd]
    &&
    Y^{[2]}
    \ar[rr]
    \ar[dd]^g
    &&
    X
    \ar[dd]
    \\
    \\
    G
    \ar[rr]
    &&
    \mathrm{INN}(G)
    \ar[rr]
    &&
    \mathbf{B} G
    \ar[rr]
    &&
    \{\bullet\}
  }
  }
  \,,
$$
produces the bundle of groupoids
$$
 \xymatrix{
    Y^{[2]} \times_{g} \mathrm{INN}(G)
    \ar[rr]
    &&
    Y^{[2]}  
  }
$$
which plays the role of the total space of the $G$-bundle classified by
$g$.
\end{observation}
This should be compared with the simplicial constructions described, for instance,
in \cite{Jurco}.

\paragraph{Remark.} Using the fact that $\mathrm{INN}(G)$ is a 2-group,
and using the injection $G \to \mathrm{INN}(G)$ we naturally
obtain the $G$-action on $Y^{[2]} \times_{g} \mathrm{INN}(G)$.

\paragraph{Remark.} Notice that this is closely related to the 
integrated Atiyah sequence 
$$
  \xymatrix{
    \mathrm{Ad}P \ar[r]
    &
    P \times_G P
    \ar[r]
    &
    X\times X
  }
$$
of groupoids over $X\times X$ coming from the $G$-principal bundle $P \to X$:
$$
  \raisebox{40pt}{
  \xymatrix{
    \mathrm{Ad}P \ar[rr]
    &&
    P \times_G P
    \ar[rr]
    &&
    X\times X
    \\
    Y \times G
    \ar[rr]
    \ar[dd]
    &&
    Y^{[2]} \times_{g} \mathrm{INN}(G)
    \ar[rr]
    \ar[dd]
    &&
    Y^{[2]}
    \ar[rr]
    \ar[dd]^g
    &&
    X
    \ar[dd]
    \\
    \\
    G
    \ar[rr]
    &&
    \mathrm{INN}(G)
    \ar[rr]
    &&
    \mathbf{B} G
    \ar[rr]
    &&
    \{\bullet\}
  }
  }
  \,.
$$

We now make precise in which sense, in turn, $Y^{[2]} \times_g \mathrm{INN}(G)$
plays the role of the total space of the $G$-bundle characterized by the
cocycle $g$.

To reobtain the $G$-bundle $P \to X$ from the groupoid 
$Y^{[2]} \times_g \mathrm{INN}(G)$ we form the pushout of
\begin{eqnarray}
  \label{first pushout}
  \xymatrix{
    Y^{[2]}\times_g \mathrm{INN}(G)
    \ar[rr]^{\mathrm{target}}
    \ar[dd]_{\mathrm{source}}
    &&
    Y \times G
    \\
    \\
    Y \times G
  }
  \,.
\end{eqnarray}

\begin{proposition}
  If $g$ is the cocycle classifying a $G$-bundle $P$ on $X$, then
  the pushout of \ref{first pushout} is (up to isomorphism)
  that very $G$-bundle $P$.
\end{proposition}
\proof
 Consider the square
 $$
  \raisebox{40pt}{
  \xymatrix{
    Y^{[2]}\times_g \mathrm{INN}(G)
    \ar[rr]^{\mathrm{target}}
    \ar[dd]_{\mathrm{source}}
    &&
    Y \times G
    \ar[d]^{t^{-1}}
    \\
    && 
    \pi^* P
    \ar[d]
    \\
    Y \times G
    \ar[r]^{t^{-1}}
    &
    \pi^* P
    \ar[r]
    &
    P
  }
  }\,,
  $$  
  where 
  $t : \xymatrix{\pi^*P \ar[r]^{\sim}&  Y \times G}$ 
  is the local trivialization of $P$ which gives rise to the
  transition function $g$. Then the diagram commutes by 
  the very definition of $g$. Since $t$ is
  an isomorphism and since
  $\pi^* P \to P$ is locally an isomorphism, it follows that this
  is the universal pushout.
\endofproof

\subsection{Universal 2-bundles in terms of $\mathrm{INN}_0(G_{(2)})$}

Now let $G_{(2)}$ be any strict 2-group. Let $Y^{[3]}$ be the
2-groupoid whose 2-morphisms are triples of lifts to $Y$ of points in $X$.
A principal $G_{(2)}$-2-bundle \cite{Bartels, Bakovic} has local 
trivializations characterized by 2-functors
$$
  g : Y^{[3]} \to \mathbf{B} G_{(2)}
  \,.
$$
\begin{definition}[$G_{(2)}$-cocycles]
  A $G_{(2)}$-(2-)cocycle on $X$ is a 2-functor
  $$
    g : Y^{[3]} \to \mathbf{B} G_{(2)}
    \,.
  $$
\end{definition}
(Instead of 2-functors on $Y^{[3]}$ one could use pseudo functors on $Y^{[2]}$.)

As before, we can pull these back along our exact sequence of 2-groupoids
\ref{the exact sequence}
$$
  \xymatrix{
    G_{(2)}
    \ar[r]
    &
    \mathrm{INN}(G_{(2)})
    \ar[r]
    &
    \mathbf{B} G_{(2)}
  }
$$
to obtain
$$
  \raisebox{40pt}{
  \xymatrix{
    Y \times G_{(2)}
    \ar[rr]
    \ar[dd]
    &&
    Y^{[3]} \times_{g} \mathrm{INN}(G_{(2)})
    \ar[rr]
    \ar[dd]
    &&
    Y^{[3]}
    \ar[rr]
    \ar[dd]^g
    &&
    X
    \ar[dd]
    \\
    \\
    G_{(2)}
    \ar[rr]
    &&
    \mathrm{INN}(G_{(2)})
    \ar[rr]
    &&
    \mathbf{B} G_{(2)}
    \ar[rr]
    &&
    \{\bullet\}
  }
  }
  \,.
$$

We reconstruct the total 2-space of the 2-bundle by forming the
weak pushout of
\begin{eqnarray}
  \label{second pushout}
  \raisebox{30pt}{
  \xymatrix{
    Y^{[3]}\times_g \mathrm{INN}(G_{(2)})
    \ar[rr]^{\mathrm{target}}
    \ar[dd]_{\mathrm{source}}
    &&
    Y \times G_{(2)}
    \\
    \\
    Y \times G_{(2)}
  }
  }
  \,.
\end{eqnarray}
Here ``source'' and ``target'' are defined relative to the inclusion
$$
  Y \times G_{(2)} \hookrightarrow Y^{[2]} \times_g \mathrm{INN}(G_{(2)})
  \,.
$$
This means that for a given 1-morphism
$$
\raisebox{36pt}
{
\xymatrix{
&&
\bullet
\ar[dd]|(.6){f := g(x,y)}
&
x
\ar[dd]
\\
\bullet 
 \ar @/^1pc/ [urr]^{q_1}_(.35){\ }="s" 
 \ar @/_1pc/ [drr]_{q_2}^(.35){\ }="t"_{}&&\\
&&\bullet 
&
y
\ar@{=>} "s";"t"^{F}
}
}
$$
in $Y^{[3]} \times_g \mathrm{INN}(G_{(2)})$ (for any $x,y \in Y$ with $\pi(x) = \pi(y)$
and for $g(x,y)$ the corresponding component of the given 2-cocycle)
which we may equivalently rewrite as
$$
\raisebox{36pt}
{
\xymatrix{
&&
\bullet
\ar[dd]^{f}
&
x
\ar[dd]
\\
\bullet 
 \ar @/^1pc/ [urr]^{q_1}_(.35){\ }="s" 
 \ar @/_1pc/ [urr]_{f^{-1} q_2}^(.35){\ }="t"_{}&&\\
&&\bullet 
&y
\ar@{=>} "s";"t"^{\act{f^{-1}}F}
}
}
$$
the source in this sense is
$$
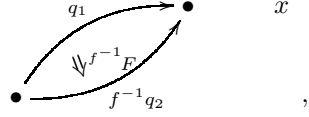
\raisebox{36pt}
{
\xymatrix{
&&
\bullet
&
x
\\
\bullet 
 \ar @/^1pc/ [urr]^{q_1}_(.35){\ }="s" 
 \ar @/_1pc/ [urr]_{f^{-1} q_2}^(.35){\ }="t"_{}&&\\
%&&\bullet 
\ar@{=>} "s";"t"^{\act{f^{-1}}F}
}
}
\,,
$$
regarded as a morphism in $Y \times G_{(2)}$,
while the target is
$$
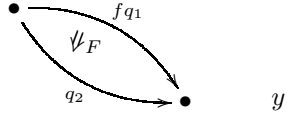
\raisebox{36pt}
{
\xymatrix{
&&
&
\\
\bullet 
 \ar @/^1pc/ [drr]^{f q_1}_(.35){\ }="s" 
 \ar @/_1pc/ [drr]_{q_2}^(.35){\ }="t"_{}&&\\
&&\bullet 
&
y
\ar@{=>} "s";"t"^{F}
}
}
$$
regarded as a morphism in $Y \times G_{(2)}$.

This way the transition function $g(x,y)$ acts on the 
copies of $G_{(2)}$ which appear as the trivialized
fibers of the $G_{(2)}$-bundle.

Bartels \cite{Bartels}[proof of prop. 22] gives a reconstruction of total space of 
principal 2-bundle from their 2-cocycles which is closely related
to $Y^{[3]} \times_g \mathrm{INN}(G_{(2)})$.

As an anonymous referee pointed out, our construction in the case of 1-groups is related 
to the universal cover. Since $\mathbf{B}G$ is a (model of a) connected 1-type, 
and $\mathrm{INN}_0(G)$ is contractible, it can be considered as (a model of) the universal 
cover. Indeed, one of the motivations for the first author was to understand 2-connected 
`universal covers' for 2-types. The connections with such a notion will be treated 
in \cite{Roberts}, as well as generalisations.

\subsection{Relation to simplicial bundles}

\label{relation to simplicial bundles}

Our considerations can be translated, along the nerve or double nerve
functor, to the world of simplicial sets. Under this translation one
finds that the tangent category construction corresponds to the simplicial
operation known as d{\'e}calage, and $\mathrm{INN}_0(G)$ corresponds to the
simplicial set denoted $W G$. It does not appear to be well known that a group structure can be put on $WG$, and one way of seeing this for simplicial groups which are nerves of 2-groups is via our construction. There is a general description of this group structure bypassing $\mathrm{INN}_0(G)$ altogether \cite{RStevenson}.

\subsubsection{Tangent categories and d{\'e}calage}

For $C$ any category, notice that a sequence of $k$ composable morphisms
$$
  
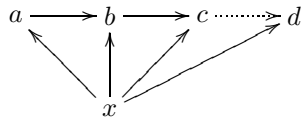
\xymatrix{
    a 
    \ar[r]
    &
    b
    \ar[r]
    &
    c
    \ar@{..>}[r]
    &
    d
    \\
    & x
    \ar[ul]
    \ar[u]
    \ar[ur]
    \ar[urr]
  }
$$
in $TC$ is, since all triangles commute, the same as a sequence 
$$
  \xymatrix{
    a 
    \ar[r]
    &
    b
    \ar[r]
    &
    c
    \ar@{..>}[r]
    &
    d
    \\
    & x
    \ar[ul]
  }
$$
of $k+1$ composable morphisms in $C$.

To formalize this observation, 
let $\Delta$ denote, as usual the simplicial category whose objects are the 
catgories
$$
  [n] = \{0 \to 1 \to \cdots \to n \}
$$
for all $n \in \mathbb{N}$
and whose morphisms are the functors between these. Write
$[\Delta^{op},\mathrm{Set}]$  for the category of simplicial sets.

  Denote by
  $$
    [(-) + 1] : \Delta^{\mathrm{op}} \to \Delta^{\mathrm{op}}
  $$
  the obvious functor which acts on objects as
  $$
    [n] \mapsto [n+1]
  $$
  (shifts everything up by one). The induced map on simplicial sets, 
  $$
    [(-)+1]^*:[\Delta^{\mathrm{op}},\mathrm{Set}] \to [\Delta^{\mathrm{op}},\mathrm{Set}]
  $$
   is called \emph{d{\'e}calage} \cite{Illusie} and is denoted $\mathrm{Dec}^1$. From a more pedestrian viewpoint, $\mathrm{Dec}^1$ strips off the first face and first 
degeneracy map\footnote{It is a matter of convention that the first face and degeneracy maps are removed. 
The d\'ecalage is sometimes defined by removing the last face and degeneracy maps, but 
our tangent category construction is related to the former convention.} from each level of a 
simplicial object $X$, re-indexes the rest and moves the sets of simplices down one level:
$$
 (\mathrm{Dec}^1 X)_n = X_{n+1}
 \,.
$$

\begin{proposition}
  The tangent category construction from definition \ref{tangent 2-bundle}
  is taken by the nerve functor $N : \mathrm{Cat} \to [\Delta^{\mathrm{op}},\mathrm{Set}]$ 
  to the d{\'e}calage construction in that
  we have a weakly commuting square
  $$
    \raisebox{40pt}{
    \xymatrix{
      \mathrm{Cat}
      \ar[dd]^T
      \ar[rr]^N
      &&
      [\Delta^{\mathrm{op}},\mathrm{Set}]
      \ar[dd]^{[(-) + 1]^*}_<{\ }="s"
      \\
      \\
      \mathrm{Cat}
      \ar[rr]_N^<{\ }="t"
      &&
      [\Delta^{\mathrm{op}},\mathrm{Set}]
      \ar@{=>}^{\simeq} "s"; "t"
    }
    }
    \,.
  $$
\end{proposition}
\proof 
  This is not hard to see by chasing explicit elements through this diagram.
  A little more abstractly, we see as follows that the assignment of $(n+1)$-simplices in $C$ to
  $n$-simplices in $TC$ is functorial.
 
  Notice that $n$-simplices in $T C$ are commuting squares
  $$
    \xymatrix{
      [n]
      \ar[rr]
      \ar[dd]
      &&
      I \times [n]
      \ar[dd]
      \\
      \\
      [0]
      \ar[rr] &&C
    }
  $$
  but the pushout of this co-cone is $[n+1]$:
  $$
    \raisebox{30pt}{
    \xymatrix{
      [n]
      \ar[rr]
      \ar[dd]
      &&
      I \times [n]
      \ar[dd]^{f}
      \\
      \\
      [0]
      \ar[rr] && [n+1]
    }
    }
    \,,
  $$
  where 
  $$
    \hspace{-2cm}
    f 
    \hspace{6pt}
      :
    \hspace{6pt}
    \raisebox{20pt}{
    \xymatrix@C=8pt{
      (\circ,0) 
      \ar[r]
      &
      (\circ,1)
      \ar[r]
      &
      (\circ,2)
      \ar@{..>}[r]
      &
      (\circ,n)
      \\
      (\bullet,0)
      \ar[r]
      \ar[u]
      &
      (\bullet,1)
      \ar[r]
      \ar[u]
      &
      (\bullet,2)
      \ar@{..>}[r]
      \ar[u]
      &
      (\bullet,n)
      \ar[u]
    }
    }
    \hspace{3pt}
      \mapsto
    \hspace{3pt}
    \raisebox{20pt}{
    \xymatrix@C=8pt{
      (\circ,0)
      \ar[r]
      &
      (\circ,1)
      \ar[r]
      &
      (\circ,2)
      \ar@{..>}[r]
      &
      (\circ,k)
      \\
      &
      (\bullet,0)
      \ar[ul]
      \ar[u]
      \ar[ur]
      \ar[urr]
    }
    }
    \hspace{3pt}
      \mapsto
    \hspace{3pt}
    \raisebox{20pt}{
    \xymatrix@C=8pt{
      (\circ,0)
      \ar[r]
      &
      (\circ,1)
      \ar[r]
      &
      (\circ,2)
      \ar@{..>}[r]
      &
      (\circ,n)
      \\
      &
      (\bullet,0)
      \ar[ul]
    }
    }
    \,.
  $$
  Hence we functorially assign $(n+1)$-simplices in $C$ to $n$-simplices in $TC$ by
  using the universality of the pushout:
  $$
    \raisebox{50pt}{
    \xymatrix{
      [n]
      \ar[rr]
      \ar[dd]
      &&
      I \times [n]
      \ar[dd]
      \ar[dl]_f
      \\
       & [n+1]
       \ar[dr]^{!}
      \\
      [0]
      \ar[rr] 
      \ar[ur]
      && C
    }
    }
    \,.
  $$
\endofproof

\subsubsection{Universal simplicial bundles}
\label{Universal simplicial bundles}

If $G$ is a simplicial group in sets, there is a notion of principal bundle internal to 
$s\mathrm{Set}$ \cite{May}, and for such bundles there is a classifying simplicial set 
$\overline{W}G$ completely analogous to the case of topological bundles. 
As such, there is a contractible simplicial set $WG$ which is the total 
space of a simplicial bundle 
$$
  \xymatrix{
     W G 
     \ar[d]
     \\
     \overline{W}G
  }
$$
and in fact $WG \simeq \mathrm{Dec}^1\overline{W}G$.

It is a short calculation to show that $\overline{W}G = N\mathbf{B}G$ 
and $WG = N\mathrm{INN}(G)$ 
when G is a constant simplicial group. To recover a similar result 
for strict 2-groups (in $\mathrm{Set}$), we 
recall that for 2-categories there is a functor $\mathcal{N}$ 
called the \emph{double nerve}, 
which forms a bisimplicial set whose geometric realization is called the classifying 
space of the 2-category.

\begin{definition}
  Recalling that strict 2-groups are the same as categories internal
  to groups, the double nerve $\mathcal{N} G_{(2)}$ of a strict 2-group 
  $G_{(2)}$ is defined to be its
  image under
  $$
    \mathcal{N}
    :
    \xymatrix{
      \mathrm{Cat}(\mathrm{Gp})
      \ar[r]^{N}
      &
      [\Delta^{\mathrm{op}}, \mathrm{Gp}]
      \ar[r]
      &
      [\Delta^{\mathrm{op}}, [\Delta^{\mathrm{op}}, \mathrm{Set}]]
    }
    \,,   
  $$
  with groups considered as one-object groupoids.
\end{definition}
Explicitly, let $G_{(2)}$ be the strict 2-group coming from the crossed module $t : H \to G$. 
Recalling that then $\mathrm{Obj}(G_{(2)}) = G$ and $\mathrm{Mor}(G_{(2)}) = H \rtimes G$ 
we find that $\mathcal{N}G_{(2)}$ is the bisimplicial set given by
\begin{eqnarray*}
  (\mathcal{N} G_{(2)})_{0n} &&= \{\cdot\}
  \\
  (\mathcal{N} G_{(2)})_{1n}
  &&=
  (N G_{(2)})_n
  =
  \left\lbrace
    \begin{array}{ll}
      G, & n = 0;
      \\
      H \rtimes G &  n = 1;
      \\
      (H \rtimes G) \times_G \stackrel{n}{\cdots} \times_G (H \rtimes G), & n > 1
    \end{array}
  \right.
  \\
  (\mathcal{N}G_{(2)})_{k n}
  &&=
  (\mathcal{N} G_{(2)})_{1n} 
    \times \stackrel{k}{\cdots} \times
  (\mathcal{N} G_{(2)})_{1n}  
  \hspace{14pt}
  k > 1
  \,.
\end{eqnarray*}
So for each $n$, $(\mathcal{N}G_{(2)})_{\bullet n}$ is the nerve of the group of 
sequences of $n$ composable morphisms in $G_{(2)}$. Applying $\mathrm{Dec}^1$ to each of
these, i.e. forming the bisimplicial set $\mathrm{Dec}^1 \mathcal{N} G_{(2)}$ in the image
of
  $$
    \mathrm{Dec}^1\mathcal{N}
    :
    \xymatrix{
      \mathrm{Cat}(\mathrm{Gp})
      \ar[r]^{N}
      &
      [\Delta^{\mathrm{op}}, \mathrm{Gp}]
      \ar[r]
      &
      [\Delta^{\mathrm{op}}, [\Delta^{\mathrm{op}}, \mathrm{Set}]]
      \ar[rr]^{[\Delta^{\mathrm{op}}, \mathrm{Dec}^1]}
      &&
      [\Delta^{\mathrm{op}}, [\Delta^{\mathrm{op}}, \mathrm{Set}]]
    }
  $$
 we obtain a surjection $\mathrm{Dec}^1 \mathcal{N} G_{(2)} \to \mathcal{N} G_{(2)}$
 whose kernel is the bisimplicial set which has just the fiber, namely the group of sequences of $n$-morphisms,
 in each row:
 $$
   (N' G_{(2)})_{k n} := (N G_{(2)})_n \hspace{20pt} \forall k \in \mathbb{N}
   \,.
 $$
 Hence we have an exact sequence
 $$
   N' G_{(2)} \to \mathrm{Dec}^1\mathcal{N}G_{(2)} \to \mathcal{N} G_{(2)}
   \,.
 $$
 The realization $|\cdot|$ of a bisimplicial set is the ordinary realization of the
diagonal simplicial space $\underline{n} \mapsto (\mathcal{N}G_{(2)})_{nn}$.  It is hence clear 
that
 $$
   |N' G_{(2)}| = |N G_{(2)}|
 $$
 and, by definition,
 $$
   |\mathcal{N} G_{(2)}| = B G_{(2)}
   \,.
 $$
  Moreover, since each row of $\mathrm{Dec}^1 \mathcal{N} G_{(2)}$ is contractible, 
  $|\mathrm{Dec}^1 \mathcal{N} G_{(2)}|$ is also contractible.\footnote{The realization of a bisimplicial set can also be calculated as first taking the ordinary realization of the rows then realizing the resulting simplicial space - this results in a space homeomorphic to the previous description.}

This relates to our construction by the fact that $\mathrm{Dec}^1 \mathcal{N} G_{(2)}$ is the nerve of a double groupoid (in fact a $cat^2$-group \cite{Loday}) constructed from the crossed square 
$$
  \raisebox{37pt}{
  \xymatrix{
   H 
       \ar[rr]^{\mathrm{id}}
       \ar[dd]_{t}
       && 
       H
       \ar[dd]^{t}
      \\
      \\
      G \ar[rr]^{\mathrm{id}}&& G 
    }
    }
    \,.
  $$

\appendix

\section{Crossed squares}

As noted above, crossed squares were introduced in \cite{Guin-Walery--Loday}. 
We include the definition for completeness:

\begin{definition}[crossed square]
  A crossed square is a commutative square of $P$-groups
$$ 
  \xymatrix{
  L \ar[rr]^f \ar[dd]_u & & M \ar[dd]^v \\
  \\
  N \ar[rr]^g & & P
  }
$$
(with $P$ acting on itself by conjugation, actions denoted by $\act{a}{b}$) and a function $M \times N \to L$ such that
\begin{enumerate}
\item $f$ and $u$ are $P$-equivariant, and 
\begin{eqnarray*}
N \longrightarrow P \\
M \longrightarrow P \\
L \longrightarrow P
\end{eqnarray*}
are crossed modules,

\item $f(h(x,y)) = x\act{g(y)}{x^{-1}}$, \hspace{11ex} $u(h(x,y)) = \act{v(x)}{y}y^{-1}$,

\item $h(f(z),y) = z\act{g(y)}{z^{-1}}$, \hspace{11ex} $h(x,u(z)) = \act{v(x)}{z}z^{-1}$,

\item $h(xx',y) = \act{v(x)}{h(x',y)}h(x,y)$, \hspace {3ex} $h(x,yy') = h(x,y)\act{g(y)}{h(x,y')}$,

\item $h$ is $P$-equivariant for the diagonal action of $P$ on $M\times N$.

\end{enumerate}
\end{definition}

It follows from the definition that $L \to M$ and $L \to N$ are crossed modules where $M,N$ act on $L$ via the maps to $P$. The reader is encouraged to note the analogies between these axioms and those in Definition \ref{2xmod definition}.

%%%%%%%%%%%%

\medskip

  \noindent {\sc School of Mathematical Sciences, The University of Adelaide, SA, 5005, Australia}
  \\
  \noindent{\it E-mail address}: {\tt droberts@maths.adelaide.edu.au}\\
  
  \noindent {\sc Department Mathematik, Universit{\"a}t Hamburg, Bundesstra{\ss}e 55,
     20146 Hamburg, Germany}\\
  \noindent {\it E-mail address}: {\tt urs.schreiber@math.uni-hamburg.de}


\begin{thebibliography}{10}

\bibitem{ArvasiUlualan} Z. Arvasi, E. Ulualan, 
\emph{On algebraic models for homotopy 3-types}, JHRS vol. 1 no. 1, pp. 1--17, (2006)

\bibitem{BaezLauda} J. Baez, A. Lauda, \emph{Higher-Dimensional Algebra V: 2-Groups},
Th. and App. of Cat. 12 (2004), pp. 423--491, available as
\verb"math/0307200".

\bibitem{BaezS} J. Baez, U. Schreiber, \emph{Higher Gauge Theory},
in \emph{Noncommutative Geometry and Representation Theory in
Mathematical Physics}, Contemporary Mathematics 391,
availabvle as \verb"math/0511710".

\bibitem{BaezStevenson} J. Baez, D. Stevenson, 
\emph{The classifying space of a topological 2-group}, available as
\verb"arXiv:0801.3843"

\bibitem{Bakovic} I. Bakovi{\'c}, \emph{Bigroupoid bitorsors}, PhD Thesis.

\bibitem{Bartels} T. Bartels, \emph{2-Bundles}, PhD Thesis, available as 
\verb"math/0410328".

\bibitem{BreenMessing} L. Breen, W. Messing, 
\emph{Differential geometry of gerbes}, Adv. in Mathematics, vol 198, 2, 
pp. 732--846 (2005)

\bibitem{BrownIcen} R. Brown, I. Icen, 
\emph{Homotopies and automorphisms of crossed modules of groupoids},
Appl. Categ. Structures, 11 (2)  (2003) pp. 185--203


\bibitem{BrownGilbert} R. Brown, N. D. Gilbert, 
\emph{Algebraic models of 3-types and automorphism structures for crossed modules},
Proc. London Math. Soc. (3) 59 (1)  (1989) pp. 51--73

\bibitem{BrownSpencer} R. Brown, C. B. Spencer, 
\emph{$\mathcal{G}$-groupoids, crossed modules and the classifying space of a topological group},
Proc. Kon. Akad. v. Wet. 79, pp. 296--302


\bibitem{Carrasco-Cegarra}P. Carrasco, A. M. Cegarra, \emph{Group-theoretic algebraic models for homotopy types}, J. Pure Appl. Algebra \textbf{75} (1991) pp. 195--235.

\bibitem{Conduche_84}D. Conduch{\'e}, \emph{Modules crois\'es g\'en\'eralis\'es de longueur 2},  J. Pure Appl. Algebra 
\textbf{34} (1984), pp. 155--178.

\bibitem{Conduche letter}
D. Conduch{\'e}, \emph{Unpublished letter to R. Brown and J.\!-L. Loday}, 1984.

\bibitem{Conduche}
D. Conduch{\'e}, \emph{Simplicial crossed modules and mapping cones}, Georgian Mathematical Journal \textbf{10} (2003) pp. 623--636.

\bibitem{Day-Street}
B.Day, R. Street, \emph{Monoidal bicategories and {H}opf algebroids}, Adv. Math \textbf{129} (1997) pp. 99--157.

\bibitem{GordonPowerStreet}
R. Gordon, A. J. Power, R. Street, 
\emph{Coherence for tricategories}  
Memoirs of the AMS, vol. 117, no. 558. (1995)

\bibitem{Gray} 
J. W. Gray, \emph{Formal Category Theory: Adjointness for 2-Categories}
Springer (1974).

\bibitem{Guin-Walery--Loday}
D. Guin-Wal\'ery and J.-L. Loday, \emph{Obstructions \`a l'Excision en K-th\'eorie 
Alg\'ebrique}, Springer Lecture Notes in Math., 854, (1981), pp. 179--216. 

\bibitem{Illusie}L. Illusie, \emph{Complexe cotangent et d\'eformation II}, Lecture Notes in Mathematics 283, Springer-Verlag, Berlin-New York, (1972).


\bibitem{Jurco}
B. Jur{\v c}o, \emph{Crossed Module Bundle Gerbes; Classification, 
String Group and Differential Geometry}, available as 
\verb"math/0510078".

\bibitem{KampsPorterBook} K. H. Kamps and T. Porter,
\emph{Abstract Homotopy and Simple Homotopy Theory}, World Scientific (1997).

\bibitem{KampsPorter} K. H. Kamps and T. Porter,
\emph{2-Groupoid enrichments in homotopy theory and algebra},
K-Theory 25, pp. 373--409 (2002).

\bibitem{Loday}
J.\!-L. Loday, \emph{Spaces with finitely many homotopy groups}, J. Pure Appl. Algebra \textbf{24} (1982) pp. 179--202.

\bibitem{May}J. P. May, \emph{Simplicial objects in algebraic topology} Van Nostrand (1967).


\bibitem{Norrie}K. Norrie, \emph{Actions and automorphisms of crossed modules.} Bulletin de la Soci\'et\'e Math\'ematique de France, \textbf{118} no. 2 (1990), pp. 129--146.

\bibitem{Roberts} D. M. Roberts, work in progress.

\bibitem{RStevenson} D. M. Roberts and D. Stevenson, \emph{On realizations of simplicial bundles},
in preparation.


\bibitem{LieStruc}
H. Sati, U. Schreiber, J. Stasheff, \emph{$L_\infty$ algebra connections and applications
to String- and Chern-Simons $n$-transport}, available as \verb"arXiv:0801.3480".

\bibitem{SWaldorfI} U. Schreiber, K. Waldorf, \emph{Parallel transport and
functors}, available as \verb"arXiv:0705.0452".

\bibitem{SWaldorfII} U. Schreiber, K. Waldorf, 
\emph{Smooth 2-Functors vs. Differential Forms}, available as \verb"arXiv:0802.0663".

\bibitem{Segal} G. B. Segal, \emph{Classifying spaces and spectral sequences},
Publ. Math. IHES No. 34 (1968) pp. 105--112

\bibitem{Whitehead} J. H. C. Whitehead, \emph{Note on a previous paper entitled
``On adding relations to homotopy groups´´, Ann. Math. {\bf 47}, pp. 806--810 (1946)}


\end{thebibliography}
\end{document}